\documentclass[a4paper]{article}
\usepackage{amsfonts}
\usepackage{amsmath}

\input{tcilatex}

\begin{document}

\title{Appendix to \\
``The Chow rings of generalized Grassmannians'' }
\author{Haibao Duan and Xuezhi Zhao \\
%EndAName
Institute of Mathematics, Chinese Academy of Sciences\\
Department of Mathematics, Capital Normal University}
\maketitle

\begin{abstract}
In this appendix we tabulate preliminary data required to establish Theorem
1-14 in the paper ``The Chow rings of generalized Grassmannians''.
\end{abstract}

\date{}

This appendix is set to record the intermediate data used in establishing
Theorem 1-14 in [DZ$_{3}$]. At the first sight, the data may be seen as
many. However, they are all produced by computer, and can be summarized in
the tables contained in the proofs of Theorem 8--14 in [DZ$_{3}$,\textbf{\ }%
\S 5], and [DZ$_{3}$,\textbf{\ }\S 8]

\bigskip

Here is a brief guide for the content tabulated. Sections \S 1--\S 7 are
devoted, respectively, to the seven cases of $(G,H)$ concerned by Theorem
1-7 (resp. Theorem 8--14) in [DZ$_{3}$].

\begin{quote}
a) In \textbf{1.1--7.1}, the sets $W(H;G)$ (cf. [DZ$_{3}$,\textbf{\ 2.1}])
of left cosets of the Weyl group of $G$ by the Weyl group of $H$ are
presented both in terms of minimized decompositions of its elements, and
index system [DZ$_{3}$,(2.3)] imposed by the decompositions. These were
generated using the \textsl{Decomposition} in [DZ$_{1}$].

b) In \textbf{1.2--7.2}, the matrices $A_{k}$ (cf. [DZ$_{3}$, \textbf{4.3}, 
\textbf{Step 1}]), $k\geq 1$, required to determine the additive cohomology
of $G/H_{s}$ are listed. They were computed using the \textsl{L--R
coefficients }in [DZ$_{1}$].

c) In \textbf{1.3--7.3}, the multiplicative rule for the basis elements of $%
H^{even}(G/H_{s})$ are presented (cf. [DZ$_{3}$, \textbf{4.3}, \textbf{Step 2%
}]). They were also obtained by applying the \textsl{L--R coefficients }in
[DZ$_{1}$].

d) In \textbf{1.4--7.4}, with respect to the ordered monomial basis $B(m)$
(cf. [DZ$_{3}$, \textbf{3.3}]) for the degrees $m$ that are relevant to the
computations in [DZ$_{3}$, \S 6], the corresponding structure matrices $%
M(\pi _{m})$ were computed by using the \textsl{L--R coefficients }in [DZ$%
_{1}$].

e) Applying the built--in function ``\textsl{Nullspace}'' in \textsl{%
Mathematica }to the $M(\pi _{m})$ in \textbf{1.4--7.4}, one obtains the 
\textsl{Nullspace} $N(\pi _{m})$ listed in \textbf{1.5--7.5}. These are
required to specify the relations $r_{i}$'s in the proofs of Theorem 1--7 in
[DZ$_{3}$,\S 6].
\end{quote}

\bigskip

We conclude this introduction with the Cartan matrices for the Lie groups of
types $F_{4}$, $E_{6}$, $E_{7}$ and $E_{8}$ [Hu, p.59]. As has been
emphasized in \S 1--\S 4 of [DZ$_{3}$], all data in this text are generated
from them.

\begin{center}
$F_{4}:\left( 
\begin{array}{cccc}
2 & -1 & 0 & 0 \\ 
-1 & 2 & -2 & 0 \\ 
0 & -1 & 2 & -1 \\ 
0 & 0 & -1 & 2%
\end{array}%
\right) $;$\qquad E_{6}:\left( 
\begin{array}{cccccc}
2 & 0 & -1 & 0 & 0 & 0 \\ 
0 & 2 & 0 & -1 & 0 & 0 \\ 
-1 & 0 & 2 & -1 & 0 & 0 \\ 
0 & -1 & -1 & 2 & -1 & 0 \\ 
0 & 0 & 0 & -1 & 2 & -1 \\ 
0 & 0 & 0 & 0 & -1 & 2%
\end{array}%
\right) $;

$E_{7}:\left( 
\begin{array}{ccccccc}
2 & 0 & -1 & 0 & 0 & 0 & 0 \\ 
0 & 2 & 0 & -1 & 0 & 0 & 0 \\ 
-1 & 0 & 2 & -1 & 0 & 0 & 0 \\ 
0 & -1 & -1 & 2 & -1 & 0 & 0 \\ 
0 & 0 & 0 & -1 & 2 & -1 & 0 \\ 
0 & 0 & 0 & 0 & -1 & 2 & -1 \\ 
0 & 0 & 0 & 0 & 0 & -1 & 2%
\end{array}%
\right) $

$E_{8}:\left( 
\begin{array}{cccccccc}
2 & 0 & -1 & 0 & 0 & 0 & 0 & 0 \\ 
0 & 2 & 0 & -1 & 0 & 0 & 0 & 0 \\ 
-1 & 0 & 2 & -1 & 0 & 0 & 0 & 0 \\ 
0 & -1 & -1 & 2 & -1 & 0 & 0 & 0 \\ 
0 & 0 & 0 & -1 & 2 & -1 & 0 & 0 \\ 
0 & 0 & 0 & 0 & -1 & 2 & -1 & 0 \\ 
0 & 0 & 0 & 0 & 0 & -1 & 2 & -1 \\ 
0 & 0 & 0 & 0 & 0 & 0 & -1 & 2%
\end{array}%
\right) $
\end{center}

\section{$(G,H_{s})=(F_{4},C_{3})$}

\textbf{1.1. The minimized decompositions of elements in }$W(C_{3},F_{4})$:

\begin{center}
\begin{tabular}{|r|l|r|}
\hline
& $w_{i,j}$ & decomposition \\ \hline
$1$ & $w_{1,1}$ & $[{1}]$ \\ \hline
$2$ & $w_{2,1}$ & $[{2, 1}]$ \\ \hline
$3$ & $w_{3,1}$ & $[{3, 2, 1}]$ \\ \hline
$4$ & $w_{4,1}$ & $[{2, 3, 2, 1}]$ \\ \hline
$5$ & $w_{4,2}$ & $[{4, 3, 2, 1}]$ \\ \hline
$6$ & $w_{5,1}$ & $[{1, 2, 3, 2, 1}]$ \\ \hline
$7$ & $w_{5,2}$ & $[{2, 4, 3, 2, 1}]$ \\ \hline
$8$ & $w_{6,1}$ & $[{1, 2, 4, 3, 2, 1}]$ \\ \hline
$9$ & $w_{6,2}$ & $[{3, 2, 4, 3, 2, 1}]$ \\ \hline
$10$ & $w_{7,1}$ & $[{1, 3, 2, 4, 3, 2, 1}]$ \\ \hline
$11$ & $w_{7,2}$ & $[{2, 3, 2, 4, 3, 2, 1}]$ \\ \hline
\end{tabular}
\begin{tabular}{|r|l|r|}
\hline
$12$ & $w_{8,1}$ & $[{1, 2, 3, 2, 4, 3, 2, 1}]$ \\ \hline
$13$ & $w_{8,2}$ & $[{2, 1, 3, 2, 4, 3, 2, 1}]$ \\ \hline
$14$ & $w_{9,1}$ & $[{1, 2, 1, 3, 2, 4, 3, 2, 1}]$ \\ \hline
$15$ & $w_{9,2}$ & $[{3, 2, 1, 3, 2, 4, 3, 2, 1}]$ \\ \hline
$16$ & $w_{10,1}$ & $[{1, 3, 2, 1, 3, 2, 4, 3, 2, 1}]$ \\ \hline
$17$ & $w_{10,2}$ & $[{4, 3, 2, 1, 3, 2, 4, 3, 2, 1}]$ \\ \hline
$18$ & $w_{11,1}$ & $[{1, 4, 3, 2, 1, 3, 2, 4, 3, 2, 1}]$ \\ \hline
$19$ & $w_{11,2}$ & $[{2, 1, 3, 2, 1, 3, 2, 4, 3, 2, 1}]$ \\ \hline
$20$ & $w_{12,1}$ & $[{2, 1, 4, 3, 2, 1, 3, 2, 4, 3, 2, 1}]$ \\ \hline
$21$ & $w_{13,1}$ & $[{3, 2, 1, 4, 3, 2, 1, 3, 2, 4, 3, 2, 1}]$ \\ \hline
$22$ & $w_{14,1}$ & $[{2, 3, 2, 1, 4, 3, 2, 1, 3, 2, 4, 3, 2, 1}]$ \\ \hline
$23$ & $w_{15,1}$ & $[{1, 2, 3, 2, 1, 4, 3, 2, 1, 3, 2, 4, 3, 2, 1}]$ \\ 
\hline
\end{tabular}
\end{center}

\textbf{1.2. Product with the Euler class }$\omega =s_{1,1}$

\begin{center}
$A_{2} = \left( 
\begin{array}{c}
1%
\end{array}%
\right) $ $A_{3} = \left( 
\begin{array}{c}
2%
\end{array}%
\right) $ $A_{4} = \left( 
\begin{array}{cc}
1 & 2%
\end{array}%
\right) $ $A_{5} = \left( 
\begin{array}{cc}
1 & 2 \\ 
0 & 1%
\end{array}%
\right) $ $A_{6} = \left( 
\begin{array}{cc}
2 & 0 \\ 
1 & 2%
\end{array}%
\right) $ $A_{7} = \left( 
\begin{array}{cc}
2 & 0 \\ 
1 & 1%
\end{array}%
\right) $ $A_{8} = \left( 
\begin{array}{cc}
1 & 2 \\ 
2 & 1%
\end{array}%
\right) $ $A_{9} = \left( 
\begin{array}{cc}
1 & 0 \\ 
1 & 2%
\end{array}%
\right) $ $A_{10} = \left( 
\begin{array}{cc}
2 & 0 \\ 
1 & 2%
\end{array}%
\right) $ $A_{11} = \left( 
\begin{array}{cc}
2 & 1 \\ 
1 & 0%
\end{array}%
\right) $ $A_{12} = \left( 
\begin{array}{c}
1 \\ 
2%
\end{array}%
\right) $ $A_{13} = \left( 
\begin{array}{c}
2%
\end{array}%
\right) $ $A_{14} = \left( 
\begin{array}{c}
1%
\end{array}%
\right) $ $A_{15} = \left( 
\begin{array}{c}
1%
\end{array}%
\right) $
\end{center}

\textbf{The additive basis of} $H^{\ast }(F_{4}/C_{3})$ derived from the $%
A_{k}$'s:

\begin{center}
\begin{tabular}{|l|l|l|}
\hline
$H^{0}$ & $Z$ &  \\ \hline
$H^{6}$ & $Z_{2}$ & $\overline{s}_{3,1}$ \\ \hline
$H^{8}$ & $Z$ & $\overline{s}_{4,2}$ \\ \hline
$H^{12}$ & $Z_{4}$ & $\overline{s}_{6,2}$ \\ \hline
$H^{14}$ & $Z_{2}$ & $\overline{s}_{7,1}$ \\ \hline
$H^{16}$ & $Z_{3}$ & $\overline{s}_{8,1}$ \\ \hline
\end{tabular}
\begin{tabular}{|l|l|l|}
\hline
$H^{18}$ & $Z_{2}$ & $\overline{s}_{9,2}$ \\ \hline
$H^{20}$ & $Z_{4}$ & $\overline{s}_{10,2}$ \\ \hline
$H^{26}$ & $Z_{2}$ & $\overline{s}_{13,1}$ \\ \hline
$H^{23}$ & $Z$ & $\beta ^{-1}(-2\,s_{11,1}+s_{11,2})$ \\ \hline
$H^{31}$ & $Z$ &  \\ \hline
&  &  \\ \hline
\end{tabular}
\end{center}

\textbf{1.3. The ring structure on} $H^{even}(F_{4}/C_{3})$:

\begin{center}
\begin{tabular}{l}
$\overline{{s}}{_{3,1}}^{2}=-2\,\overline{s}_{6,2}$, \\ 
$\overline{s}_{3,1}\,\overline{s}_{4,2}=\overline{s}_{7,1}$, \\ 
$\overline{s}_{3,1}\,\overline{s}_{6,2}=3\,\overline{s}_{9,2}$,%
\end{tabular}
\ 
\begin{tabular}{l}
$\overline{s}_{3,1}\,\overline{s}_{7,1}=-6\,\overline{s}_{10,2}$, \\ 
$\overline{s}_{10,2}\,\overline{s}_{3,1}=\overline{s}_{13,1}$, \\ 
$\overline{{s}}{_{4,2}}^{2}=-\overline{s}_{8,1}$,%
\end{tabular}%
\ 
\begin{tabular}{l}
$\overline{s}_{4,2}\,\overline{s}_{6,2}=-3\,\overline{s}_{10,2}$, \\ 
$\overline{s}_{4,2}\,\overline{s}_{9,2}=\overline{s}_{13,1}$, \\ 
$\overline{s}_{6,2}\,\overline{s}_{7,1}=3\,\overline{s}_{13,1}$.%
\end{tabular}
\end{center}

\textbf{1.4. The structure matrices }$M(\pi _{m})$\textbf{\ }are with
respect to the ordered monomial base $B(2m)$ (cf. [DZ, 3.3])

$%
\begin{array}{rcl}
B(6) & = & \{y_{3},{y_{1}^{3}}\}\text{;}%
\end{array}%
$

$%
\begin{array}{rcl}
B(12) & = & \{y_{6},y_{3}^{2},{y_{1}^{2}}y_{4},{y_{1}^{3}}y_{3},{y_{1}^{6}}\}%
\end{array}%
$

$%
\begin{array}{rcl}
B(16) & = & \{{y_{4}},y_{1}y_{3}y_{4},{y_{1}^{2}}y_{6},\ {y_{1}^{2}y_{3}^{2}}%
,{y_{1}^{4}}y_{4},\ {y_{1}^{5}}y_{3},{y_{1}^{8}}\}%
\end{array}%
$

$%
\begin{array}{rcl}
B(24) & = & \{{y}_{6}^{2},\ {y_{4}^{3}},\ {y_{3}^{2}}y_{6},\ \ {y_{3}^{4}},\
y_{1}y_{3}{y_{4}^{2}},\ \ {y_{1}^{2}}y_{4}y_{6},\ \ {y_{1}^{2}}\,{y_{3}^{2}}%
y_{4},\ \ {y_{1}^{3}}y_{3}y_{6},\ {y_{1}^{3}y_{3}^{3}}, \\ 
&  & \ \ {y_{1}^{4}y_{4}^{2}},\ {y_{1}^{5}}y_{3,1}y_{4},\ \ {y_{1}^{6}}%
y_{6},\ {y_{1}^{6}y_{3}^{2}},\ \ {y_{1}^{8}}y_{4},\ {y_{1}^{9}}\,y_{3},\ \ {%
y_{1}^{12}}\}%
\end{array}%
$

\noindent of\textsl{\ }$\mathbb{Z}[y_{1},y_{3},y_{4},y_{6}]^{(2m)}$ for $%
m=3,6,8,12$ respectively, where $%
y_{1}=s_{1,1},y_{3}=s_{3,1},y_{4}=s_{4,2},y_{6}=s_{6,2}$.

\begin{center}
$M(\pi _{3})=\left( 
\begin{array}{cc}
1 & 2%
\end{array}%
\right) $;$\quad M(\pi _{6})=\left( 
\begin{array}{ccccc}
0 & 3 & 1 & 6 & 12 \\ 
1 & 4 & 2 & 8 & 16%
\end{array}%
\right) $;

$M(\pi _{8})=\left( 
\begin{array}{ccccccc}
1 & 5 & 3 & 24 & 10 & 48 & 96 \\ 
1 & 4 & 3 & 18 & 8 & 36 & 72%
\end{array}%
\right) $;

{\small \setlength{\arraycolsep}{1.2mm} $M(\pi _{12})=(%
\begin{array}{cccccccccccccccc}
3 & 3 & 21 & 156 & 14 & 9 & 66 & 42 & 312 & 28 & 132 & 84 & 624 & 264 & \
1248 & 2496%
\end{array}%
)$. }
\end{center}

\textbf{1.5. The Nullspace of }$M(\pi _{m})$ for $m=3,6,8,12$ are
respectively

\begin{center}
$N(\pi _{3})=\left( 
\begin{array}{cc}
-2 & 1%
\end{array}%
\right) $; $\ N(\pi _{6})=\left( 
\begin{array}{ccccc}
0 & -4 & 0 & 0 & 1 \\ 
0 & -2 & 0 & 1 & 0 \\ 
-2 & -1 & 3 & 0 & 0%
\end{array}%
\right) $;

$N(\pi _{8})=\left( 
\begin{array}{ccccccc}
24 & -24 & 0 & 0 & 0 & 0 & 1 \\ 
12 & -12 & 0 & 0 & 0 & 1 & 0 \\ 
0 & -2 & 0 & 0 & 1 & 0 & 0 \\ 
6 & -6 & 0 & 1 & 0 & 0 & 0 \\ 
-3 & 0 & 1 & 0 & 0 & 0 & 0%
\end{array}%
\right) $;

$N(\pi _{12})=\left( 
\begin{array}{cccccccccccccccc}
-832 & 0 & 0 & 0 & 0 & 0 & 0 & 0 & 0 & 0 & 0 & 0 & 0 & 0 & 0 & 1 \\ 
-416 & 0 & 0 & 0 & 0 & 0 & 0 & 0 & 0 & 0 & 0 & 0 & 0 & 0 & 1 & 0 \\ 
-88 & 0 & 0 & 0 & 0 & 0 & 0 & 0 & 0 & 0 & 0 & 0 & 0 & 1 & 0 & 0 \\ 
-208 & 0 & 0 & 0 & 0 & 0 & 0 & 0 & 0 & 0 & 0 & 0 & 1 & 0 & 0 & 0 \\ 
-28 & 0 & 0 & 0 & 0 & 0 & 0 & 0 & 0 & 0 & 0 & 1 & 0 & 0 & 0 & 0 \\ 
-44 & 0 & 0 & 0 & 0 & 0 & 0 & 0 & 0 & 0 & 1 & 0 & 0 & 0 & 0 & 0 \\ 
-28 & 0 & 0 & 0 & 0 & 0 & 0 & 0 & 0 & 3 & 0 & 0 & 0 & 0 & 0 & 0 \\ 
-104 & 0 & 0 & 0 & 0 & 0 & 0 & 0 & 1 & 0 & 0 & 0 & 0 & 0 & 0 & 0 \\ 
-14 & 0 & 0 & 0 & 0 & 0 & 0 & 1 & 0 & 0 & 0 & 0 & 0 & 0 & 0 & 0 \\ 
-22 & 0 & 0 & 0 & 0 & 0 & 1 & 0 & 0 & 0 & 0 & 0 & 0 & 0 & 0 & 0 \\ 
-3 & 0 & 0 & 0 & 0 & 1 & 0 & 0 & 0 & 0 & 0 & 0 & 0 & 0 & 0 & 0 \\ 
-14 & 0 & 0 & 0 & 3 & 0 & 0 & 0 & 0 & 0 & 0 & 0 & 0 & 0 & 0 & 0 \\ 
-52 & 0 & 0 & 1 & 0 & 0 & 0 & 0 & 0 & 0 & 0 & 0 & 0 & 0 & 0 & 0 \\ 
-7 & 0 & 1 & 0 & 0 & 0 & 0 & 0 & 0 & 0 & 0 & 0 & 0 & 0 & 0 & 0 \\ 
-1 & 1 & 0 & 0 & 0 & 0 & 0 & 0 & 0 & 0 & 0 & 0 & 0 & 0 & 0 & 0%
\end{array}%
\right) $.
\end{center}

\section{$(G,H_{s})=(F_{4},B_{3})$}

\textbf{2.1. The minimized decompositions of elements in }$W(B_{3},F_{4})$:

\begin{center}
\begin{tabular}{|r|l|r|}
\hline
& $w_{i,j}$ & decomposition \\ \hline
$1$ & $w_{1,1}$ & $[{4}]$ \\ \hline
$2$ & $w_{2,1}$ & $[{3, 4}]$ \\ \hline
$3$ & $w_{3,1}$ & $[{2, 3, 4}]$ \\ \hline
$4$ & $w_{4,1}$ & $[{1, 2, 3, 4}]$ \\ \hline
$5$ & $w_{4,2}$ & $[{3, 2, 3, 4}]$ \\ \hline
$6$ & $w_{5,1}$ & $[{1, 3, 2, 3, 4}]$ \\ \hline
$7$ & $w_{5,2}$ & $[{4, 3, 2, 3, 4}]$ \\ \hline
$8$ & $w_{6,1}$ & $[{1, 4, 3, 2, 3, 4}]$ \\ \hline
$9$ & $w_{6,2}$ & $[{2, 1, 3, 2, 3, 4}]$ \\ \hline
$10$ & $w_{7,1}$ & $[{2, 1, 4, 3, 2, 3, 4}]$ \\ \hline
$11$ & $w_{7,2}$ & $[{3, 2, 1, 3, 2, 3, 4}]$ \\ \hline
\end{tabular}
\begin{tabular}{|r|l|r|}
\hline
$12$ & $w_{8,1}$ & $[{3, 2, 1, 4, 3, 2, 3, 4}]$ \\ \hline
$13$ & $w_{8,2}$ & $[{4, 3, 2, 1, 3, 2, 3, 4}]$ \\ \hline
$14$ & $w_{9,1}$ & $[{2, 3, 2, 1, 4, 3, 2, 3, 4}]$ \\ \hline
$15$ & $w_{9,2}$ & $[{3, 4, 3, 2, 1, 3, 2, 3, 4}]$ \\ \hline
$16$ & $w_{10,1}$ & $[{1, 2, 3, 2, 1, 4, 3, 2, 3, 4}]$ \\ \hline
$17$ & $w_{10,2}$ & $[{2, 3, 4, 3, 2, 1, 3, 2, 3, 4}]$ \\ \hline
$18$ & $w_{11,1}$ & $[{1, 2, 3, 4, 3, 2, 1, 3, 2, 3, 4}]$ \\ \hline
$19$ & $w_{11,2}$ & $[{3, 2, 3, 4, 3, 2, 1, 3, 2, 3, 4}]$ \\ \hline
$20$ & $w_{12,1}$ & $[{1, 3, 2, 3, 4, 3, 2, 1, 3, 2, 3, 4}]$ \\ \hline
$21$ & $w_{13,1}$ & $[{2, 1, 3, 2, 3, 4, 3, 2, 1, 3, 2, 3, 4}]$ \\ \hline
$22$ & $w_{14,1}$ & $[{3, 2, 1, 3, 2, 3, 4, 3, 2, 1, 3, 2, 3, 4}]$ \\ \hline
$23$ & $w_{15,1}$ & $[{4, 3, 2, 1, 3, 2, 3, 4, 3, 2, 1, 3, 2, 3, 4}]$ \\ 
\hline
\end{tabular}
\end{center}

\textbf{2.2. Product with the Euler class }$\omega =s_{1,1}$

\begin{center}
$A_{2} = \left( 
\begin{array}{c}
1%
\end{array}%
\right) $ $A_{3} = \left( 
\begin{array}{c}
1%
\end{array}%
\right) $ $A_{4} = \left( 
\begin{array}{cc}
1 & 1%
\end{array}%
\right) $ $A_{5} = \left( 
\begin{array}{cc}
1 & 0 \\ 
1 & 1%
\end{array}%
\right) $ $A_{6} = \left( 
\begin{array}{cc}
1 & 1 \\ 
1 & 0%
\end{array}%
\right) $ $A_{7} = \left( 
\begin{array}{cc}
1 & 0 \\ 
1 & 1%
\end{array}%
\right) $ $A_{8} = \left( 
\begin{array}{cc}
2 & 1 \\ 
1 & 2%
\end{array}%
\right) $ $A_{9} = \left( 
\begin{array}{cc}
1 & 1 \\ 
0 & 1%
\end{array}%
\right) $ $A_{10} = \left( 
\begin{array}{cc}
1 & 1 \\ 
0 & 1%
\end{array}%
\right) $ $A_{11} = \left( 
\begin{array}{cc}
1 & 0 \\ 
1 & 1%
\end{array}%
\right) $ $A_{12} = \left( 
\begin{array}{c}
1 \\ 
1%
\end{array}%
\right) $ $A_{13} = \left( 
\begin{array}{c}
1%
\end{array}%
\right) $ $A_{14} = \left( 
\begin{array}{c}
1%
\end{array}%
\right) $ $A_{15} = \left( 
\begin{array}{c}
1%
\end{array}%
\right) $
\end{center}

\textbf{The additive basis of }$H^{\ast }(F_{4}/B_{3})$ derived from the $%
A_{k}$'s:

\begin{center}
\begin{tabular}{|l|l|l|l|}
\hline
$H^{0}$ & $Z$ &  &  \\ \hline
$H^{8}$ & $Z$ & $\bar{s}_{4,2}$ &  \\ \hline
$H^{16}$ & $Z_{3}$ & $\bar{s}_{8,1}$ &  \\ \hline
$H^{23}$ & $Z$ & $\beta ^{-1}(-s_{11,1}+s_{11,2})$ &  \\ \hline
$H^{31}$ & $Z$ &  &  \\ \hline
\end{tabular}
\end{center}

\textbf{2.3. The ring structure on }$H^{even}(F_{4}/B_{3})$\textbf{. }

\begin{center}
${\bar{s}_{4,2}}^{2}=\bar{s}_{8,1}$.
\end{center}

\textbf{2.4. The structure matrices }$M(\pi _{m})$\textbf{\ }are with
respect to the ordered monomial base

\begin{center}
$%
\begin{array}{rcl}
B(16) & = & \{{y}_{4}^{2},\ {y}_{1}^{4}\,y_{4},\ {y_{1}^{8}}\}\text{;}%
\end{array}%
\begin{array}{rcl}
B(24) & = & \{{y_{4}^{3}},\ {y}_{1}^{4}\,{y}_{4}^{2},\ {y_{1}^{8}}\,y_{4},\ {%
y_{1}^{12}}\}%
\end{array}%
$
\end{center}

\noindent respectively (cf. [DZ$_{3}$, \textbf{3.3}]), where ${y_{1}=s}%
_{1,1} $, $\,y_{4}=s_{4,2}$,

\begin{center}
\begin{equation*}
M(\pi_{8}) = \left( 
\begin{array}{ccc}
3 & 5 & 9 \\ 
4 & 7 & 12%
\end{array}%
\right) \ \ M(\pi_{12}) = \left( 
\begin{array}{cccc}
15 & 26 & 45 & 78%
\end{array}%
\right)
\end{equation*}
\end{center}

\textbf{2.5. The Nullspaces of }$M(\pi _{m})$ for $m=8,12$ are respectively

\begin{equation*}
N(\pi_{8}) = \left( 
\begin{array}{ccc}
-3 & 0 & 1%
\end{array}%
\right) \ \ N(\pi_{12}) = \left( 
\begin{array}{cccc}
-26 & 0 & 0 & 5 \\ 
-3 & 0 & 1 & 0 \\ 
-26 & 15 & 0 & 0%
\end{array}%
\right)
\end{equation*}

\section{$(G,H_{s})=(E_{6},A_{6})$}

\textbf{3.1. The minimized decompositions of elements in }$W(A_{6};E_{6})$:

\begin{center}
{\small \setlength{\tabcolsep}{1mm} \lineskip 0mm }
\end{center}

{\small 
\begin{tabular}{|r|l|r|}
\hline
& $w_{i,j}$ & decomposition \\ \hline
$1$ & $w_{1,1}$ & $[{2}]$ \\ \hline
$2$ & $w_{2,1}$ & $[{4, 2}]$ \\ \hline
$3$ & $w_{3,1}$ & $[{3, 4, 2}]$ \\ \hline
$4$ & $w_{3,2}$ & $[{5, 4, 2}]$ \\ \hline
$5$ & $w_{4,1}$ & $[{1, 3, 4, 2}]$ \\ \hline
$6$ & $w_{4,2}$ & $[{3, 5, 4, 2}]$ \\ \hline
$7$ & $w_{4,3}$ & $[{6, 5, 4, 2}]$ \\ \hline
$8$ & $w_{5,1}$ & $[{1, 3, 5, 4, 2}]$ \\ \hline
$9$ & $w_{5,2}$ & $[{3, 6, 5, 4, 2}]$ \\ \hline
$10$ & $w_{5,3}$ & $[{4, 3, 5, 4, 2}]$ \\ \hline
$11$ & $w_{6,1}$ & $[{1, 3, 6, 5, 4, 2}]$ \\ \hline
$12$ & $w_{6,2}$ & $[{1, 4, 3, 5, 4, 2}]$ \\ \hline
$13$ & $w_{6,3}$ & $[{2, 4, 3, 5, 4, 2}]$ \\ \hline
$14$ & $w_{6,4}$ & $[{4, 3, 6, 5, 4, 2}]$ \\ \hline
$15$ & $w_{7,1}$ & $[{1, 2, 4, 3, 5, 4, 2}]$ \\ \hline
$16$ & $w_{7,2}$ & $[{1, 4, 3, 6, 5, 4, 2}]$ \\ \hline
$17$ & $w_{7,3}$ & $[{2, 4, 3, 6, 5, 4, 2}]$ \\ \hline
$18$ & $w_{7,4}$ & $[{3, 1, 4, 3, 5, 4, 2}]$ \\ \hline
$19$ & $w_{7,5}$ & $[{5, 4, 3, 6, 5, 4, 2}]$ \\ \hline
$20$ & $w_{8,1}$ & $[{1, 2, 4, 3, 6, 5, 4, 2}]$ \\ \hline
$21$ & $w_{8,2}$ & $[{1, 5, 4, 3, 6, 5, 4, 2}]$ \\ \hline
$22$ & $w_{8,3}$ & $[{2, 3, 1, 4, 3, 5, 4, 2}]$ \\ \hline
\end{tabular}
\ 
\begin{tabular}{|r|l|r|}
\hline
$23$ & $w_{8,4}$ & $[{2, 5, 4, 3, 6, 5, 4, 2}]$ \\ \hline
$24$ & $w_{8,5}$ & $[{3, 1, 4, 3, 6, 5, 4, 2}]$ \\ \hline
$25$ & $w_{9,1}$ & $[{1, 2, 5, 4, 3, 6, 5, 4, 2}]$ \\ \hline
$26$ & $w_{9,2}$ & $[{2, 3, 1, 4, 3, 6, 5, 4, 2}]$ \\ \hline
$27$ & $w_{9,3}$ & $[{3, 1, 5, 4, 3, 6, 5, 4, 2}]$ \\ \hline
$28$ & $w_{9,4}$ & $[{4, 2, 3, 1, 4, 3, 5, 4, 2}]$ \\ \hline
$29$ & $w_{9,5}$ & $[{4, 2, 5, 4, 3, 6, 5, 4, 2}]$ \\ \hline
$30$ & $w_{10,1}$ & $[{1, 4, 2, 5, 4, 3, 6, 5, 4, 2}]$ \\ \hline
$31$ & $w_{10,2}$ & $[{2, 3, 1, 5, 4, 3, 6, 5, 4, 2}]$ \\ \hline
$32$ & $w_{10,3}$ & $[{3, 4, 2, 5, 4, 3, 6, 5, 4, 2}]$ \\ \hline
$33$ & $w_{10,4}$ & $[{4, 2, 3, 1, 4, 3, 6, 5, 4, 2}]$ \\ \hline
$34$ & $w_{10,5}$ & $[{4, 3, 1, 5, 4, 3, 6, 5, 4, 2}]$ \\ \hline
$35$ & $w_{10,6}$ & $[{5, 4, 2, 3, 1, 4, 3, 5, 4, 2}]$ \\ \hline
$36$ & $w_{11,1}$ & $[{1, 3, 4, 2, 5, 4, 3, 6, 5, 4, 2}]$ \\ \hline
$37$ & $w_{11,2}$ & $[{2, 4, 3, 1, 5, 4, 3, 6, 5, 4, 2}]$ \\ \hline
$38$ & $w_{11,3}$ & $[{3, 1, 4, 2, 5, 4, 3, 6, 5, 4, 2}]$ \\ \hline
$39$ & $w_{11,4}$ & $[{4, 2, 3, 1, 5, 4, 3, 6, 5, 4, 2}]$ \\ \hline
$40$ & $w_{11,5}$ & $[{5, 4, 2, 3, 1, 4, 3, 6, 5, 4, 2}]$ \\ \hline
$41$ & $w_{11,6}$ & $[{6, 5, 4, 2, 3, 1, 4, 3, 5, 4, 2}]$ \\ \hline
$42$ & $w_{12,1}$ & $[{1, 3, 1, 4, 2, 5, 4, 3, 6, 5, 4, 2}]$ \\ \hline
$43$ & $w_{12,2}$ & $[{2, 4, 2, 3, 1, 5, 4, 3, 6, 5, 4, 2}]$ \\ \hline
$44$ & $w_{12,3}$ & $[{3, 4, 2, 3, 1, 5, 4, 3, 6, 5, 4, 2}]$ \\ \hline
$45$ & $w_{12,4}$ & $[{4, 5, 4, 2, 3, 1, 4, 3, 6, 5, 4, 2}]$ \\ \hline
\end{tabular}
}

{\small 
\begin{tabular}{|r|l|r|}
\hline
$46$ & $w_{12,5}$ & $[{5, 6, 5, 4, 2, 3, 1, 4, 3, 5, 4, 2}]$ \\ \hline
$47$ & $w_{13,1}$ & $[{1, 3, 4, 2, 3, 1, 5, 4, 3, 6, 5, 4, 2}]$ \\ \hline
$48$ & $w_{13,2}$ & $[{2, 3, 4, 2, 3, 1, 5, 4, 3, 6, 5, 4, 2}]$ \\ \hline
$49$ & $w_{13,3}$ & $[{2, 4, 5, 4, 2, 3, 1, 4, 3, 6, 5, 4, 2}]$ \\ \hline
$50$ & $w_{13,4}$ & $[{3, 4, 5, 4, 2, 3, 1, 4, 3, 6, 5, 4, 2}]$ \\ \hline
$51$ & $w_{13,5}$ & $[{4, 5, 6, 5, 4, 2, 3, 1, 4, 3, 5, 4, 2}]$ \\ \hline
$52$ & $w_{14,1}$ & $[{1, 2, 3, 4, 2, 3, 1, 5, 4, 3, 6, 5, 4, 2}]$ \\ \hline
$53$ & $w_{14,2}$ & $[{1, 3, 4, 5, 4, 2, 3, 1, 4, 3, 6, 5, 4, 2}]$ \\
\hline
\end{tabular}
}

\textbf{3.2. Product with the Euler class }$\omega =s_{1,1}$

\begin{center}
$A_{2}=\left( 
\begin{array}{c}
1%
\end{array}%
\right) $, $A_{3}=\left( 
\begin{array}{cc}
1 & 1%
\end{array}%
\right) $, $A_{4}=\left( 
\begin{array}{ccc}
1 & 1 & 0 \\ 
0 & 1 & 1%
\end{array}%
\right) $, $A_{5}=\left( 
\begin{array}{ccc}
1 & 0 & 0 \\ 
1 & 1 & 1 \\ 
0 & 1 & 0%
\end{array}%
\right) $,

$A_{6}=\left( 
\begin{array}{cccc}
1 & 1 & 0 & 0 \\ 
1 & 0 & 0 & 1 \\ 
0 & 1 & 1 & 1%
\end{array}%
\right) $, $A_{7}=\left( 
\begin{array}{ccccc}
0 & 1 & 0 & 0 & 0 \\ 
1 & 1 & 0 & 1 & 0 \\ 
1 & 0 & 1 & 0 & 0 \\ 
0 & 1 & 1 & 0 & 1%
\end{array}%
\right) $, $A_{8}=\left( 
\begin{array}{ccccc}
1 & 0 & 1 & 0 & 0 \\ 
1 & 1 & 0 & 0 & 1 \\ 
1 & 0 & 0 & 1 & 0 \\ 
0 & 0 & 1 & 0 & 1 \\ 
0 & 1 & 0 & 1 & 0%
\end{array}%
\right) $,$A_{9}=\left( 
\begin{array}{ccccc}
1 & 1 & 0 & 0 & 0 \\ 
1 & 0 & 1 & 0 & 0 \\ 
0 & 1 & 0 & 1 & 0 \\ 
1 & 0 & 0 & 0 & 1 \\ 
0 & 1 & 1 & 0 & 0%
\end{array}%
\right) $, $A_{10}=\left( 
\begin{array}{cccccc}
1 & 1 & 0 & 0 & 0 & 0 \\ 
0 & 1 & 0 & 1 & 0 & 0 \\ 
0 & 1 & 0 & 0 & 1 & 0 \\ 
0 & 0 & 0 & 1 & 0 & 1 \\ 
1 & 0 & 1 & 0 & 0 & 0%
\end{array}%
\right) $, \ $A_{11}=\left( 
\begin{array}{cccccc}
1 & 0 & 2 & 1 & 0 & 0 \\ 
0 & 1 & 1 & 2 & 1 & 0 \\ 
2 & 0 & 1 & 0 & 0 & 0 \\ 
0 & 0 & 0 & 1 & 2 & 1 \\ 
0 & 2 & 0 & 1 & 0 & 0 \\ 
0 & 0 & 0 & 0 & 1 & 2%
\end{array}%
\right) $,

$A_{12}= \left( 
\begin{array}{ccccc}
1 & 0 & 0 & 0 & 0 \\ 
0 & 1 & 0 & 0 & 0 \\ 
1 & 0 & 1 & 0 & 0 \\ 
0 & 1 & 1 & 1 & 0 \\ 
0 & 0 & 0 & 1 & 1 \\ 
0 & 0 & 0 & 0 & 1%
\end{array}%
\right) $, \ $A_{13}= \left( 
\begin{array}{ccccc}
1 & 0 & 0 & 0 & 0 \\ 
0 & 1 & 1 & 0 & 0 \\ 
1 & 1 & 0 & 1 & 0 \\ 
0 & 0 & 1 & 1 & 1 \\ 
0 & 0 & 0 & 0 & 1%
\end{array}%
\right) $, \ $A_{14}= \left( 
\begin{array}{ccccc}
1 & 1 & 0 & 0 & 0 \\ 
1 & 0 & 1 & 0 & 0 \\ 
0 & 0 & 1 & 1 & 0 \\ 
0 & 1 & 1 & 0 & 1 \\ 
0 & 0 & 0 & 1 & 1%
\end{array}%
\right) $,

$A_{15}= \left( 
\begin{array}{cccc}
1 & 0 & 0 & 0 \\ 
1 & 1 & 0 & 0 \\ 
1 & 0 & 1 & 1 \\ 
0 & 0 & 1 & 0 \\ 
0 & 1 & 1 & 0%
\end{array}%
\right) $, \ $A_{16}= \left( 
\begin{array}{ccc}
1 & 1 & 0 \\ 
1 & 0 & 0 \\ 
1 & 0 & 1 \\ 
0 & 1 & 1%
\end{array}%
\right) $, \ $A_{17}= \left( 
\begin{array}{ccc}
1 & 0 & 0 \\ 
1 & 1 & 0 \\ 
1 & 0 & 1%
\end{array}%
\right) $,

$A_{18}= \left( 
\begin{array}{cc}
1 & 1 \\ 
0 & 1 \\ 
1 & 0%
\end{array}%
\right) $, \ $A_{19}= \left( 
\begin{array}{c}
1 \\ 
1%
\end{array}%
\right) $, \ $A_{20}= \left( 
\begin{array}{c}
1%
\end{array}%
\right) $, \ $A_{21}= \left( 
\begin{array}{c}
1%
\end{array}%
\right) $,
\end{center}

\textbf{The additive basis of} $H^{\ast }(E_{6}/A_{6})$ derived from the $%
A_{k}$'s:

\begin{center}
\begin{tabular}{|l|l|l|}
\hline
$H^{0}$ & $Z$ &  \\ \hline
$H^{6}$ & $Z$ & $\overline{s}_{3,2}$ \\ \hline
$H^{8}$ & $Z$ & $\overline{s}_{4,3}$ \\ \hline
$H^{12}$ & $Z$ & $\overline{s}_{6,1}$ \\ \hline
$H^{14}$ & $Z$ & $\overline{s}_{7,1}$ \\ \hline
\end{tabular}
\ \ \ 
\begin{tabular}{|l|l|l|}
\hline
$H^{16}$ & $Z_{3}$ & $\overline{s}_{8,1}$ \\ \hline
$H^{18}$ & $Z_{2}$ & $\overline{s}_{9,1}$ \\ \hline
$H^{20}$ & $Z$ & $\overline{s}_{10,1}$ \\ \hline
$H^{22}$ & $Z_{3}$ & $\overline{s}_{11,1}$ \\ \hline
$H^{26}$ & $Z_{2}$ & $\overline{s}_{13,2}$ \\ \hline
\end{tabular}

\begin{tabular}{|l|l|l|}
\hline
$H^{28}$ & $Z_{3}$ & $\overline{s}_{14,1}$ \\ \hline
$H^{23}$ & $Z$ & $\beta ^{-1}(s_{11,1}-s_{11,2}-s_{11,3}+s_{11,4}-\
s_{11,5}+s_{11,6})$ \\ \hline
$H^{29}$ & $Z$ & $\beta ^{-1}(-s_{14,1}+s_{14,2}+s_{14,4}-\ s_{14,5})$ \\ 
\hline
$H^{31}$ & $Z$ & $\beta ^{-1}(s_{15,1}-2\,s_{15,2}+s_{15,3}-\ s_{15,4})$ \\ 
\hline
$H^{35}$ & $Z$ & $\beta ^{-1}(-s_{17,1}+s_{17,2}+s_{17,3})$ \\ \hline
$H^{37}$ & $Z$ & $\beta ^{-1}(-s_{18,1}+s_{18,2})$ \\ \hline
$H^{43}$ & $Z$ &  \\ \hline
\end{tabular}
\end{center}

\textbf{3.3. The ring structure on} $H^{even}(E_{6}/A_{6})$:

\begin{center}
\begin{tabular}{l}
$\overline{{s}}{_{3,2}}^{2}=-2\,\overline{s}_{6,1}$ \\ 
$\overline{s}_{3,2}\,\overline{s}_{4,3}=\overline{s}_{7,1}$ \\ 
$\overline{s}_{3,2}\,\overline{s}_{6,1}=\overline{s}_{9,1}$ \\ 
$\overline{s}_{3,2}\,\overline{s}_{7,1}=2\,\overline{s}_{10,1}$ \\ 
$\overline{s}_{3,2}\,\overline{s}_{8,1}=4\,\overline{s}_{11,1}$%
\end{tabular}
\ 
\begin{tabular}{l}
$\overline{s}_{10,1}\,\overline{s}_{3,2}=-3\,\overline{s}_{13,2}$ \\ 
$\overline{s}_{11,1}\,\overline{s}_{3,2}=\overline{s}_{14,1}$ \\ 
$\overline{{s}}{_{4,3}}^{2}=\overline{s}_{8,1}$ \\ 
$\overline{s}_{4,3}\,\overline{s}_{6,1}=-\overline{s}_{10,1}$ \\ 
$\overline{s}_{4,3}\,\overline{s}_{7,1}=\,\overline{s}_{11,1}$%
\end{tabular}
\ 
\begin{tabular}{l}
$\overline{s}_{4,3}\,\overline{s}_{9,1}=\overline{s}_{13,2}$ \\ 
$\overline{s}_{10,1}\,\overline{s}_{4,3}=-\overline{s}_{14,1}$ \\ 
$\overline{s}_{6,1}\,\overline{s}_{7,1}=\overline{s}_{13,2}$ \\ 
$\overline{s}_{6,1}\,\overline{s}_{8,1}=-2\,\overline{s}_{14,1}$ \\ 
$\overline{{s}}{_{7,1}}^{2}=-2\,\overline{s}_{14,1}$%
\end{tabular}
\end{center}

\textbf{3.4. The structure matrices }$M(\pi _{m})$\textbf{\ }are with
respect to the ordered monomial base $B(2m)$

%% all monomials in the dimensions with new relations
$B(12) = \{ {y_{6}}, {y_{3}}^{2}, \ {y_{1}}^{2}{y_{4}}, {y_{1}}^{3}{y_{3}},
\ {y_{1}}^{6} \} $

$B(16) = \{ {y_{4}}^{2}, {y_{1}}{y_{3}}{y_{4}}, \ {y_{1}}^{2}{y_{6}}, {y_{1}}%
^{2}{y_{3}}^{2}, \ {y_{1}}^{4}{y_{4}}, {y_{1}}^{5}{y_{3}}, \ {y_{1}}^{8} \} $

$B(18) = \{ {y_{3}}{y_{6}}, {y_{3}}^{3}, \ {y_{1}}{y_{4}}^{2}, {y_{1}}^{2}{%
y_{3}}{y_{4}}, \ {y_{1}}^{3}{y_{6}}, {y_{1}}^{3}{y_{3}}^{2}, \ {y_{1}}^{5}{%
y_{4}}, {y_{1}}^{6}{y_{3}}, \ {y_{1}}^{9} \} $

$B(24) = \{ {y_{6}}^{2}, {y_{4}}^{3}, \ {y_{3}}^{2}{y_{6}}, {y_{3}}^{4}, \ {%
y_{1}}{y_{3}}{y_{4}}^{2}, \ {y_{1}}^{2}{y_{4}}{y_{6}}, \ {y_{1}}^{2}{y_{3}}%
^{2}{y_{4}}, \ {y_{1}}^{3}{y_{3}}{y_{6}}, \ {y_{1}}^{3}{y_{3}}^{3}$,

$\ \ \ \ \ \ \ \ \ \ \ {y_{1}}^{4}{y_{4}}^{2}, \ {y_{1}}^{5}{y_{3}}{y_{4}}, {%
y_{1}}^{6}{y_{6}}, \ {y_{1}}^{6}{y_{3}}^{2}, {y_{1}}^{8}{y_{4}}, \ {y_{1}}%
^{9}{y_{3}}, {y_{1}}^{12} \} $

\noindent for $m=6,8,9,12$, where $y_{1}={s_{1,1},y}%
_{3}=s_{3,2},y_{4}=s_{4,3},y_{6}=s_{6,1}$.

\begin{center}
$M(\pi_{6}) = \left( 
\begin{array}{ccccc}
1 & 1 & 1 & 3 & 6 \\ 
0 & 1 & 0 & 2 & 5 \\ 
0 & 0 & 0 & 1 & 2 \\ 
0 & 2 & 1 & 3 & 5 \\ 
&  &  &  & 
\end{array}%
\right) $ $M(\pi_{8}) = \left( 
\begin{array}{ccccccc}
0 & 1 & 1 & 7 & 3 & 15 & 30 \\ 
1 & 2 & 1 & 6 & 3 & 11 & 21 \\ 
0 & 0 & 0 & 2 & 0 & 5 & 12 \\ 
0 & 1 & 0 & 4 & 2 & 7 & 12 \\ 
0 & 1 & 1 & 5 & 2 & 10 & 21 \\ 
&  &  &  &  &  & 
\end{array}%
\right) $

$M(\pi_{9}) = \left( 
\begin{array}{ccccccccc}
1 & 9 & 1 & 4 & 2 & 17 & 8 & 33 & 63 \\ 
1 & 6 & 0 & 2 & 2 & 14 & 5 & 30 & 63 \\ 
1 & 6 & 1 & 3 & 2 & 11 & 5 & 21 & 42 \\ 
0 & 1 & 0 & 0 & 0 & 2 & 0 & 5 & 12 \\ 
0 & 2 & 0 & 1 & 0 & 4 & 2 & 7 & 12 \\ 
&  &  &  &  &  &  &  & 
\end{array}%
\right) $

$M(\pi_{12}) = \left( 
\begin{array}{cccccccccccccccc}
0 & 0 & 3 & 30 & 2 & 1 & 13 & 6 & 60 & 5 & 27 & 12 & 117 & 54 & 225 & 429 \\ 
1 & 1 & 7 & 50 & 5 & 3 & 22 & 14 & 99 & 10 & 43 & 28 & 196 & 84 & 390 & \ 780
\\ 
1 & 1 & 7 & 56 & 5 & 3 & 24 & 14 & 111 & 10 & 48 & 28 & 220 & 96 & 435 & \
858 \\ 
1 & 1 & 7 & 50 & 4 & 3 & 20 & 14 & 102 & 8 & 41 & 28 & 208 & 84 & 423 & \ 858
\\ 
0 & 0 & 3 & 21 & 1 & 1 & 7 & 6 & 45 & 2 & 15 & 12 & 96 & 33 & 204 & 429 \\ 
&  &  &  &  &  &  &  &  &  &  &  &  &  &  & 
\end{array}%
\right) $
\end{center}

\textbf{3.5. The Nullspace of }$M(\pi _{m})$ for $m=6,8,9,12$ are
respectively

\begin{center}
$N(\pi_{6}) = \left( 
\begin{array}{ccccc}
-2 & -1 & 3 & -2 & 1 \\ 
&  &  &  & 
\end{array}%
\right) $

$N(\pi_{8}) = \left( 
\begin{array}{ccccccc}
-3 & 6 & -3 & -6 & 3 & 0 & 1 \\ 
-3 & 6 & -1 & -5 & 0 & 2 & 0 \\ 
&  &  &  &  &  & 
\end{array}%
\right) $

$N(\pi_{9}) = \left( 
\begin{array}{ccccccccc}
-18 & -12 & 15 & 6 & 0 & 0 & 3 & 0 & 1 \\ 
-6 & -5 & 6 & 3 & 0 & 0 & 0 & 1 & 0 \\ 
-2 & -2 & 3 & 0 & 0 & 1 & 0 & 0 & 0 \\ 
-2 & 0 & 0 & 0 & 1 & 0 & 0 & 0 & 0 \\ 
&  &  &  &  &  &  &  & 
\end{array}%
\right) $

$N(\pi_{12}) = \left( 
\begin{array}{cccccccccccccccc}
26 & 0 & -52 & -26 & 0 & 0 & 39 & 0 & 0 & 0 & 0 & 0 & 0 & 0 & 0 & 1 \\ 
11 & 0 & -35 & -24 & 12 & 0 & 27 & 0 & 0 & 0 & 0 & 0 & 0 & 0 & 2 & 0 \\ 
-5 & 0 & 1 & -1 & 6 & 0 & -3 & 0 & 0 & 0 & 0 & 0 & 0 & 1 & 0 & 0 \\ 
0 & 0 & -6 & -5 & 6 & 0 & 3 & 0 & 0 & 0 & 0 & 0 & 1 & 0 & 0 & 0 \\ 
0 & 0 & -4 & 0 & 0 & 0 & 0 & 0 & 0 & 0 & 0 & 1 & 0 & 0 & 0 & 0 \\ 
-11 & 0 & -1 & 0 & 18 & 0 & -15 & 0 & 0 & 0 & 6 & 0 & 0 & 0 & 0 & 0 \\ 
0 & 0 & -2 & 1 & 0 & 0 & -3 & 0 & 0 & 3 & 0 & 0 & 0 & 0 & 0 & 0 \\ 
0 & 0 & -2 & -2 & 3 & 0 & 0 & 0 & 1 & 0 & 0 & 0 & 0 & 0 & 0 & 0 \\ 
0 & 0 & -2 & 0 & 0 & 0 & 0 & 1 & 0 & 0 & 0 & 0 & 0 & 0 & 0 & 0 \\ 
-2 & 0 & -1 & 0 & 0 & 3 & 0 & 0 & 0 & 0 & 0 & 0 & 0 & 0 & 0 & 0 \\ 
-1 & 1 & 0 & 0 & 0 & 0 & 0 & 0 & 0 & 0 & 0 & 0 & 0 & 0 & 0 & 0 \\ 
&  &  &  &  &  &  &  &  &  &  &  &  &  &  & 
\end{array}%
\right) $
\end{center}

\section{$(G,H_{s})=(E_{6},D_{5})$}

\textbf{4.1. The minimized decompositions of elements in }$W(D_{5};E_{6})$

\begin{center}
{\setlength{\tabcolsep}{1 mm} \lineskip 0mm 
\begin{tabular}{|r|l|r|}
\hline
& $w_{i,j}$ & decomposition \\ \hline
$1$ & $w_{1,1}$ & $[{6}]$ \\ \hline
$2$ & $w_{2,1}$ & $[{5,6}]$ \\ \hline
$3$ & $w_{3,1}$ & $[{4,5,6}]$ \\ \hline
$4$ & $w_{4,1}$ & $[{2,4,5,6}]$ \\ \hline
$5$ & $w_{4,2}$ & $[{3,4,5,6}]$ \\ \hline
$6$ & $w_{5,1}$ & $[{1,3,4,5,6}]$ \\ \hline
$7$ & $w_{5,2}$ & $[{2,3,4,5,6}]$ \\ \hline
$8$ & $w_{6,1}$ & $[{1,2,3,4,5,6}]$ \\ \hline
$9$ & $w_{6,2}$ & $[{4,2,3,4,5,6}]$ \\ \hline
$10$ & $w_{7,1}$ & $[{1,4,2,3,4,5,6}]$ \\ \hline
$11$ & $w_{7,2}$ & $[{5,4,2,3,4,5,6}]$ \\ \hline
$12$ & $w_{8,1}$ & $[{1,5,4,2,3,4,5,6}]$ \\ \hline
$13$ & $w_{8,2}$ & $[{3,1,4,2,3,4,5,6}]$ \\ \hline
\end{tabular}%
\begin{tabular}{|r|l|r|}
\hline
$14$ & $w_{8,3}$ & $[{6,5,4,2,3,4,5,6}]$ \\ \hline
$15$ & $w_{9,1}$ & $[{1,6,5,4,2,3,4,5,6}]$ \\ \hline
$16$ & $w_{9,2}$ & $[{3,1,5,4,2,3,4,5,6}]$ \\ \hline
$17$ & $w_{10,1}$ & $[{3,1,6,5,4,2,3,4,5,6}]$ \\ \hline
$18$ & $w_{10,2}$ & $[{4,3,1,5,4,2,3,4,5,6}]$ \\ \hline
$19$ & $w_{11,1}$ & $[{2,4,3,1,5,4,2,3,4,5,6}]$ \\ \hline
$20$ & $w_{11,2}$ & $[{4,3,1,6,5,4,2,3,4,5,6}]$ \\ \hline
$21$ & $w_{12,1}$ & $[{2,4,3,1,6,5,4,2,3,4,5,6}]$ \\ \hline
$22$ & $w_{12,2}$ & $[{5,4,3,1,6,5,4,2,3,4,5,6}]$ \\ \hline
$23$ & $w_{13,1}$ & $[{2,5,4,3,1,6,5,4,2,3,4,5,6}]$ \\ \hline
$24$ & $w_{14,1}$ & $[{4,2,5,4,3,1,6,5,4,2,3,4,5,6}]$ \\ \hline
$25$ & $w_{15,1}$ & $[{3,4,2,5,4,3,1,6,5,4,2,3,4,5,6}]$ \\ \hline
$26$ & $w_{16,1}$ & $[{1,3,4,2,5,4,3,1,6,5,4,2,3,4,5,6}]$ \\ \hline
&  &  \\ \hline
\end{tabular}
}
\end{center}

\textbf{4.2. Product with the Euler class }$\omega =s_{1,1}$

\begin{center}
$A_{2}= \left( 
\begin{array}{c}
1%
\end{array}%
\right) $, \ $A_{3}= \left( 
\begin{array}{c}
1%
\end{array}%
\right) $, \ $A_{4}= \left( 
\begin{array}{cc}
1 & 1%
\end{array}%
\right) $,

$A_{5}= \left( 
\begin{array}{cc}
0 & 1 \\ 
1 & 1%
\end{array}%
\right) $, \ $A_{6}= \left( 
\begin{array}{cc}
1 & 0 \\ 
1 & 1%
\end{array}%
\right) $, \ $A_{7}= \left( 
\begin{array}{cc}
1 & 0 \\ 
1 & 1%
\end{array}%
\right) $, \ $A_{8}= \left( 
\begin{array}{ccc}
1 & 1 & 0 \\ 
1 & 0 & 1%
\end{array}%
\right) $,

$A_{9}= \left( 
\begin{array}{cc}
1 & 1 \\ 
0 & 1 \\ 
1 & 0%
\end{array}%
\right) $, \ $A_{10}= \left( 
\begin{array}{cc}
1 & 0 \\ 
1 & 1%
\end{array}%
\right) $, \ $A_{11}= \left( 
\begin{array}{cc}
0 & 1 \\ 
1 & 1%
\end{array}%
\right) $, \ $A_{12}= \left( 
\begin{array}{cc}
1 & 0 \\ 
1 & 1%
\end{array}%
\right) $, \ $A_{13}= \left( 
\begin{array}{c}
1 \\ 
1%
\end{array}%
\right) $,

$A_{14}= \left( 
\begin{array}{c}
1%
\end{array}%
\right) $, \ $A_{15}= \left( 
\begin{array}{c}
1%
\end{array}%
\right) $, \ $A_{16}= \left( 
\begin{array}{c}
1%
\end{array}%
\right) $,
\end{center}

\textbf{The additive basis of} $H^{\ast }(E_{6}/D_{5})$ derived from the $%
A_{k}$'s:

\begin{center}
\begin{tabular}{|l|l|l|}
\hline
$H^{0}$ & $Z$ &  \\ \hline
$H^{8}$ & $Z$ & $\overline{s}_{4,1}$ \\ \hline
$H^{16}$ & $Z$ & $\overline{s}_{8,1}$ \\ \hline
\end{tabular}
\begin{tabular}{|l|l|l|}
\hline
$H^{17}$ & $Z$ & $\beta ^{-1}(-s_{8,1}+s_{8,2}+s_{8,3})$ \\ \hline
$H^{25}$ & $Z$ & $\beta ^{-1}(-s_{12,1}+s_{12,2})$ \\ \hline
$H^{33}$ & $Z$ &  \\ \hline
\end{tabular}
\end{center}

\textbf{4.3. The ring structure on} $H^{even}(E_{6}/D_{5})$:

\begin{center}
$\overline{{s}}{_{4,1}}^{2}=-\overline{s}_{8,1}$ .
\end{center}

\textbf{4.4. The structure matrices }$M(\pi _{m})$\textbf{\ }are with
respect to the ordered monomial base $B(2m)$ $:$

\begin{center}
$B(18)=\{y_{1}{y_{4}^{2}},\ {y_{1}^{5}}y_{4},\ \ {y_{1}^{9}}\}$; $\qquad
B(24)=\{{y_{4}^{3}},\ {y_{1}^{4}y_{4}^{2}},\ \ {y_{1}^{8}}y_{4},\ {y_{1}^{12}%
}\}$,
\end{center}

\noindent where $y_{1}=s_{1,1},y_{4}=s_{4,1}$:

\begin{center}
$M(\pi _{9})=\left( 
\begin{array}{ccc}
2 & 4 & 9 \\ 
2 & 5 & 12%
\end{array}%
\right) $;$\qquad M(\pi _{12})=\left( 
\begin{array}{cccc}
3 & 8 & 19 & 45 \\ 
3 & 6 & 14 & 33%
\end{array}%
\right) $.
\end{center}

\textbf{4.5. The Nullspace of }$M(\pi _{m})$ for $m=9,12$ are respectively

\begin{center}
$N(\pi _{9})=\left( 
\begin{array}{ccc}
3 & -6 & 2%
\end{array}%
\right) $;$\qquad N(\pi _{12})=\left( 
\begin{array}{cccc}
1 & -6 & 0 & 1 \\ 
2 & -15 & 6 & 0%
\end{array}%
\right) $.
\end{center}

\section{$(G,H_{s})=(E_{7},E_{6})$}

\textbf{5.1. The minimized decompositions of elements in }$W(E_{6};E_{7})$

\begin{center}
{\small \setlength{\tabcolsep}{0.3 mm} \lineskip 0mm 
\begin{tabular}{|r|l|r|}
\hline
& $w_{i,j}$ & decomposition \\ \hline
$1$ & $w_{1,1}$ & $[{7}]$ \\ \hline
$2$ & $w_{2,1}$ & $[{6,7}]$ \\ \hline
$3$ & $w_{3,1}$ & $[{5,6,7}]$ \\ \hline
$4$ & $w_{4,1}$ & $[{4,5,6,7}]$ \\ \hline
$5$ & $w_{5,1}$ & $[{2,4,5,6,7}]$ \\ \hline
$6$ & $w_{5,2}$ & $[{3,4,5,6,7}]$ \\ \hline
$7$ & $w_{6,1}$ & $[{1,3,4,5,6,7}]$ \\ \hline
$8$ & $w_{6,2}$ & $[{2,3,4,5,6,7}]$ \\ \hline
$9$ & $w_{7,1}$ & $[{1,2,3,4,5,6,7}]$ \\ \hline
$10$ & $w_{7,2}$ & $[{4,2,3,4,5,6,7}]$ \\ \hline
$11$ & $w_{8,1}$ & $[{1,4,2,3,4,5,6,7}]$ \\ \hline
$12$ & $w_{8,2}$ & $[{5,4,2,3,4,5,6,7}]$ \\ \hline
$13$ & $w_{9,1}$ & $[{1,5,4,2,3,4,5,6,7}]$ \\ \hline
$14$ & $w_{9,2}$ & $[{3,1,4,2,3,4,5,6,7}]$ \\ \hline
$15$ & $w_{9,3}$ & $[{6,5,4,2,3,4,5,6,7}]$ \\ \hline
$16$ & $w_{10,1}$ & $[{1,6,5,4,2,3,4,5,6,7}]$ \\ \hline
$17$ & $w_{10,2}$ & $[{3,1,5,4,2,3,4,5,6,7}]$ \\ \hline
$18$ & $w_{10,3}$ & $[{7,6,5,4,2,3,4,5,6,7}]$ \\ \hline
$19$ & $w_{11,1}$ & $[{1}, 7,6,5,4,2,3,4,5,6,7 ]$ \\ \hline
$20$ & $w_{11,2}$ & $[{3,1,6,5,4,2,3,4,5,6,7}]$ \\ \hline
$21$ & $w_{11,3}$ & $[{4,3,1,5,4,2,3,4,5,6,7}]$ \\ \hline
$22$ & $w_{12,1}$ & $[{2,4,3,1,5,4,2,3,4,5,6,7}]$ \\ \hline
$23$ & $w_{12,2}$ & $[{3,1}, 7,6,5,4,2,3,4,5,6,7 ]$ \\ \hline
\end{tabular}%
\begin{tabular}{|r|l|r|}
\hline
$24$ & $w_{12,3}$ & $[{4,3,1,6,5,4,2,3,4,5,6,7}]$ \\ \hline
$25$ & $w_{13,1}$ & $[{2,4,3,1,6,5,4,2,3,4,5,6,7}]$ \\ \hline
$26$ & $w_{13,2}$ & $[{4,3,1}, 7,6,5,4,2,3,4,5,6,7 ]$ \\ \hline
$27$ & $w_{13,3}$ & $[{5,4,3,1,6,5,4,2,3,4,5,6,7}]$ \\ \hline
$28$ & $w_{14,1}$ & $[{2,4,3,1}, 7,6,5,4,2,3,4,5,6,7 ]$ \\ \hline
$29$ & $w_{14,2}$ & $[{2,5,4,3,1,6,5,4,2,3,4,5,6,7}]$ \\ \hline
$30$ & $w_{14,3}$ & $[{5,4,3,1}, 7,6,5,4,2,3,4,5,6,7 ]$ \\ \hline
$31$ & $w_{15,1}$ & $[{2,5,4,3,1}, 7,6,5,4,2,3,4,5,6,7 ]$ \\ \hline
$32$ & $w_{15,2}$ & $[{4,2,5,4,3,1,6,5,4,2,3,4,5,6,7}]$ \\ \hline
$33$ & $w_{15,3}$ & $[{6,5,4,3,1}, 7,6,5,4,2,3,4,5,6,7 ]$ \\ \hline
$34$ & $w_{16,1}$ & $[{2,6,5,4,3,1}, 7,6,5,4,2,3,4,5,6,7 ]$ \\ \hline
$35$ & $w_{16,2}$ & $[{3,4,2,5,4,3,1,6,5,4,2,3,4,5,6,7}]$ \\ \hline
$36$ & $w_{16,3}$ & $[{4,2,5,4,3,1}, 7,6,5,4,2,3,4,5,6,7 ]$ \\ \hline
$37$ & $w_{17,1}$ & $[{1,3,4,2,5,4,3,1,6,5,4,2,3,4,5,6,7}]$ \\ \hline
$38$ & $w_{17,2}$ & $[{3,4,2,5,4,3,1}, 7,6,5,4,2,3,4,5,6,7 ]$ \\ \hline
$39$ & $w_{17,3}$ & $[{4,2,6,5,4,3,1}, 7,6,5,4,2,3,4,5,6,7 ]$ \\ \hline
$40$ & $w_{18,1}$ & $[{1,3,4,2,5,4,3,1}, 7,6,5,4,2,3,4,5,6,7 ]$ \\ \hline
$41$ & $w_{18,2}$ & $[{3,4,2,6,5,4,3,1}, 7,6,5,4,2,3,4,5,6,7 ]$ \\ \hline
$42$ & $w_{18,3}$ & $[{5,4,2,6,5,4,3,1}, 7,6,5,4,2,3,4,5,6,7 ]$ \\ \hline
$43$ & $w_{19,1}$ & $[{1,3,4,2,6,5,4,3,1}, 7,6,5,4,2,3,4,5,6,7 ]$ \\ \hline
$44$ & $w_{19,2}$ & $[{3,5,4,2,6,5,4,3,1}, 7,6,5,4,2,3,4,5,6,7 ]$ \\ \hline
$45$ & $w_{20,1}$ & $[{1,3,5,4,2,6,5,4,3,1}, 7,6,5,4,2,3,4,5,6,7 ]$ \\ \hline
$46$ & $w_{20,2}$ & $[{4,3,5,4,2,6,5,4,3,1}, 7,6,5,4,2,3,4,5,6,7 ]$ \\ \hline
$47$ & $w_{21,1}$ & {\ $[{1\!,4,3\!,5,4\!,2,6,5,4,3,1}, 7,6,5,4,2,3,4,5,6,7
] $} \\ \hline
\end{tabular}%
}
\end{center}

{\small $\ $ }

{\small 
\begin{tabular}{|r|l|r|}
\hline
$48$ & $w_{21,2}$ & $[{2,4,3,5,4,2,6,5,4,3,1}, 7,6,5,4,2,3,4,5,6,7 ]$ \\ 
\hline
$49$ & $w_{22,1}$ & $[{1,2,4,3,5,4,2,6,5,4,3,1}, 7,6,5,4,2,3,4,5,6,7 ]$ \\ 
\hline
$50$ & $w_{22,2}$ & $[{3,1,4,3,5,4,2,6,5,4,3,1}, 7,6,5,4,2,3,4,5,6,7 ]$ \\ 
\hline
$51$ & $w_{23,1}$ & $[{2,3,1,4,3,5,4,2,6,5,4,3,1}, 7,6,5,4,2,3,4,5,6,7 ]$ \\ 
\hline
$52$ & $w_{24,1}$ & $[{4,2,3,1,4,3,5,4,2,6,5,4,3,1}, 7,6,5,4,2,3,4,5,6,7 ]$
\\ \hline
$53$ & $w_{25,1}$ & $[{5,4,2,3,1,4,3,5,4,2,6,5,4,3,1}, 7,6,5,4,2,3,4,5,6,7 ]$
\\ \hline
$54$ & $w_{26,1}$ & $[{6,5,4,2,3,1,4,3,5,4,2,6,5,4,3,1}, 7,6,5,4,2,3,4,5,6,7
]$ \\ \hline
$55$ & $w_{27,1}$ & $[{7,6,5,4,2,3,1,4,3,5,4,2,6,5,4,3,1},
7,6,5,4,2,3,4,5,6,7 ]$ \\ \hline
\end{tabular}
}

\textbf{5.2. Product with the Euler class }$\omega =s_{1,1}$

\begin{center}
$A_{2}=\left( 
\begin{array}{c}
1%
\end{array}%
\right) $, \ $A_{3}=\left( 
\begin{array}{c}
1%
\end{array}%
\right) $, \ $A_{4}=\left( 
\begin{array}{c}
1%
\end{array}%
\right) $, \ $A_{5}=\left( 
\begin{array}{cc}
1 & 1%
\end{array}%
\right) $,

$A_{6}= \left( 
\begin{array}{cc}
0 & 1 \\ 
1 & 1%
\end{array}%
\right) $, \ $A_{7}= \left( 
\begin{array}{cc}
1 & 0 \\ 
1 & 1%
\end{array}%
\right) $, \ $A_{8}= \left( 
\begin{array}{cc}
1 & 0 \\ 
1 & 1%
\end{array}%
\right) $, \ $A_{9}= \left( 
\begin{array}{ccc}
1 & 1 & 0 \\ 
1 & 0 & 1%
\end{array}%
\right) $,

$A_{10}= \left( 
\begin{array}{ccc}
1 & 1 & 0 \\ 
0 & 1 & 0 \\ 
1 & 0 & 1%
\end{array}%
\right) $, \ $A_{11}= \left( 
\begin{array}{ccc}
1 & 1 & 0 \\ 
0 & 1 & 1 \\ 
1 & 0 & 0%
\end{array}%
\right) $, \ $A_{12}= \left( 
\begin{array}{ccc}
0 & 1 & 0 \\ 
0 & 1 & 1 \\ 
1 & 0 & 1%
\end{array}%
\right) $, \ $A_{13}= \left( 
\begin{array}{ccc}
1 & 0 & 0 \\ 
0 & 1 & 0 \\ 
1 & 1 & 1%
\end{array}%
\right) $,

$A_{14}= \left( 
\begin{array}{ccc}
1 & 1 & 0 \\ 
1 & 0 & 1 \\ 
0 & 1 & 1%
\end{array}%
\right) $, \ $A_{15}= \left( 
\begin{array}{ccc}
1 & 0 & 0 \\ 
1 & 1 & 0 \\ 
1 & 0 & 1%
\end{array}%
\right) $, \ $A_{16}= \left( 
\begin{array}{ccc}
1 & 0 & 1 \\ 
0 & 1 & 1 \\ 
1 & 0 & 0%
\end{array}%
\right) $, \ $A_{17}= \left( 
\begin{array}{ccc}
0 & 0 & 1 \\ 
1 & 1 & 0 \\ 
0 & 1 & 1%
\end{array}%
\right) $,

$A_{18}= \left( 
\begin{array}{ccc}
1 & 0 & 0 \\ 
1 & 1 & 0 \\ 
0 & 1 & 1%
\end{array}%
\right) $, \ $A_{19}= \left( 
\begin{array}{cc}
1 & 0 \\ 
1 & 1 \\ 
0 & 1%
\end{array}%
\right) $, \ $A_{20}= \left( 
\begin{array}{cc}
1 & 0 \\ 
1 & 1%
\end{array}%
\right) $, \ $A_{21}= \left( 
\begin{array}{cc}
1 & 0 \\ 
1 & 1%
\end{array}%
\right) $,

$A_{22}=\left( 
\begin{array}{cc}
1 & 1 \\ 
1 & 0%
\end{array}%
\right) $, \ $A_{23}=\left( 
\begin{array}{c}
1 \\ 
1%
\end{array}%
\right) $, \ $A_{24}=\left( 
\begin{array}{c}
1%
\end{array}%
\right) $, \ $A_{25}=\left( 
\begin{array}{c}
1%
\end{array}%
\right) $, \ $A_{26}=\left( 
\begin{array}{c}
1%
\end{array}%
\right) $, \ $A_{27}=\left( 
\begin{array}{c}
1%
\end{array}%
\right) $.
\end{center}

\textbf{The additive basis of} $H^{\ast }(E_{7}/E_{6})$ derived from the $%
A_{k}$'s:

\begin{center}
\begin{tabular}{|l|l|l|}
\hline
$H^{0}$ & $Z$ &  \\ \hline
$H^{10}$ & $Z$ & $\overline{s}_{5,1}$ \\ \hline
$H^{18}$ & $Z$ & $\overline{s}_{9,1}$ \\ \hline
$H^{28}$ & $Z_{2}$ & $\overline{s}_{14,1}$ \\ \hline
$H^{37}$ & $Z$ & $\beta ^{-1}(-s_{18,1}+s_{18,2}-s_{18,3})$ \\ \hline
$H^{45}$ & $Z$ & $\beta ^{-1}(-s_{22,1}+s_{22,2})$ \\ \hline
$H^{55}$ & $Z$ &  \\ \hline
\end{tabular}
\end{center}

\textbf{5.3. The ring structure on} $H^{even}(E_{7}/E_{6})$:

\begin{center}
$\overline{s}_{5,1}\,\overline{s}_{9,1}=3\,\overline{s}_{14,1}$ .
\end{center}

\textbf{5.4. The structure matrices }$M(\pi _{m})$\textbf{\ }are with
respect to the ordered monomial base $B(2m)$ $:$

$B(20)=\{{y_{5}^{2}},\ y_{1}y_{9},\ {y_{1}^{5}}y_{5},\ \ {y_{1}^{10}}\}$;

$B(28)=\{y_{5}y_{9},\ {y_{1}^{4}y_{5}^{2}},\ \ {y_{1}^{5}}y_{9},\ {y_{1}^{9}}%
y_{5},\ \ {y_{1}^{14}}\}$;

$B(36)=\{{y_{9}^{2}},\ {y_{1}^{3}y_{5}^{3}},\ \ {y_{1}^{4}}y_{5}y_{9},\ {%
y_{1}^{8}y_{5}^{2}},\ \ {y_{1}^{9}}y_{9},\ {y_{1}^{13}}y_{5},\ \ {y_{1}^{18}}%
\}$

\noindent where $y_{1}=s_{1,1},y_{5}=s_{5,1},y_{9}=s_{9,1}$:

\begin{center}
$M(\pi _{10})=\left( 
\begin{array}{cccc}
2 & 1 & 4 & 9 \\ 
2 & 1 & 5 & 12 \\ 
0 & 0 & 1 & 2%
\end{array}%
\right) $;$\qquad M(\pi _{14})=\left( 
\begin{array}{ccccc}
4 & 20 & 10 & 47 & 110 \\ 
3 & 14 & 7 & 33 & 78 \\ 
4 & 18 & 9 & 42 & 98%
\end{array}%
\right) $;

$M(\pi _{18})=\left( 
\begin{array}{ccccccc}
4 & 40 & 20 & 94 & 47 & 221 & 520 \\ 
10 & 92 & 46 & 216 & 108 & 507 & 1190 \\ 
6 & 58 & 29 & 136 & 68 & 319 & 748%
\end{array}%
\right) $.
\end{center}

\textbf{5.5. The Nullspace of }$M(\pi _{m})$ for $m=10,14,18$ are
respectively

\begin{center}
$N(\pi _{10})=\left( 
\begin{array}{cccc}
-1 & 2 & 0 & 0%
\end{array}%
\right) $;$\qquad N(\pi _{14})=\left( 
\begin{array}{ccccc}
-2 & 9 & 0 & -6 & 1 \\ 
0 & -1 & 2 & 0 & 0%
\end{array}%
\right) $;

$N(\pi _{18})=\left( 
\begin{array}{ccccccc}
3 & 29 & 0 & -18 & 0 & 0 & 1 \\ 
1 & 10 & 0 & -9 & 0 & 2 & 0 \\ 
0 & 0 & 0 & -1 & 2 & 0 & 0 \\ 
0 & -1 & 2 & 0 & 0 & 0 & 0%
\end{array}%
\right) $.
\end{center}

\section{$(G,H_{s})=(E_{7},D_{6})$}

\textbf{6.1. The minimized decompositions of elements in }$W(D_{6};E_{7})$

\begin{center}
{\small \setlength{\tabcolsep}{1.5mm} \lineskip 1 mm 
\begin{tabular}{|r|l|r|}
\hline
& $w_{i,j}$ & decomposition \\ \hline
$1$ & $w_{1,1}$ & $[{1}]$ \\ \hline
$2$ & $w_{2,1}$ & $[{3,1}]$ \\ \hline
$3$ & $w_{3,1}$ & $[{4,3,1}]$ \\ \hline
$4$ & $w_{4,1}$ & $[{2,4,3,1}]$ \\ \hline
$5$ & $w_{4,2}$ & $[{5,4,3,1}]$ \\ \hline
$6$ & $w_{5,1}$ & $[{2,5,4,3,1}]$ \\ \hline
$7$ & $w_{5,2}$ & $[{6,5,4,3,1}]$ \\ \hline
$8$ & $w_{6,1}$ & $[{2,6,5,4,3,1}]$ \\ \hline
$9$ & $w_{6,2}$ & $[{4,2,5,4,3,1}]$ \\ \hline
$10$ & $w_{6,3}$ & $[{7,6,5,4,3,1}]$ \\ \hline
$11$ & $w_{7,1}$ & $[{2,7,6,5,4,3,1}]$ \\ \hline
$12$ & $w_{7,2}$ & $[{3,4,2,5,4,3,1}]$ \\ \hline
$13$ & $w_{7,3}$ & $[{4,2,6,5,4,3,1}]$ \\ \hline
$14$ & $w_{8,1}$ & $[{1,3,4,2,5,4,3,1}]$ \\ \hline
$15$ & $w_{8,2}$ & $[{3,4,2,6,5,4,3,1}]$ \\ \hline
$16$ & $w_{8,3}$ & $[{4,2,7,6,5,4,3,1}]$ \\ \hline
$17$ & $w_{8,4}$ & $[{5,4,2,6,5,4,3,1}]$ \\ \hline
$18$ & $w_{9,1}$ & $[{1,3,4,2,6,5,4,3,1}]$ \\ \hline
$19$ & $w_{9,2}$ & $[{3,4,2,7,6,5,4,3,1}]$ \\ \hline
$20$ & $w_{9,3}$ & $[{3,5,4,2,6,5,4,3,1}]$ \\ \hline
$21$ & $w_{9,4}$ & $[{5,4,2,7,6,5,4,3,1}]$ \\ \hline
$22$ & $w_{10,1}$ & $[{1,3,4,2,7,6,5,4,3,1}]$ \\ \hline
$23$ & $w_{10,2}$ & $[{1,3,5,4,2,6,5,4,3,1}]$ \\ \hline
$24$ & $w_{10,3}$ & $[{3,5,4,2,7,6,5,4,3,1}]$ \\ \hline
$25$ & $w_{10,4}$ & $[{4,3,5,4,2,6,5,4,3,1}]$ \\ \hline
$26$ & $w_{10,5}$ & $[{6,5,4,2,7,6,5,4,3,1}]$ \\ \hline
\end{tabular}
\begin{tabular}{|r|l|r|}
\hline
$27$ & $w_{11,1}$ & $[{1,3,5,4,2,7,6,5,4,3,1}]$ \\ \hline
$28$ & $w_{11,2}$ & $[{1,4,3,5,4,2,6,5,4,3,1}]$ \\ \hline
$29$ & $w_{11,3}$ & $[{2,4,3,5,4,2,6,5,4,3,1}]$ \\ \hline
$30$ & $w_{11,4}$ & $[{3,6,5,4,2,7,6,5,4,3,1}]$ \\ \hline
$31$ & $w_{11,5}$ & $[{4,3,5,4,2,7,6,5,4,3,1}]$ \\ \hline
$32$ & $w_{12,1}$ & $[{1,2,4,3,5,4,2,6,5,4,3,1}]$ \\ \hline
$33$ & $w_{12,2}$ & $[{1,3,6,5,4,2,7,6,5,4,3,1}]$ \\ \hline
$34$ & $w_{12,3}$ & $[{1,4,3,5,4,2,7,6,5,4,3,1}]$ \\ \hline
$35$ & $w_{12,4}$ & $[{2,4,3,5,4,2,7,6,5,4,3,1}]$ \\ \hline
$36$ & $w_{12,5}$ & $[{3,1,4,3,5,4,2,6,5,4,3,1}]$ \\ \hline
$37$ & $w_{12,6}$ & $[{4,3,6,5,4,2,7,6,5,4,3,1}]$ \\ \hline
$38$ & $w_{13,1}$ & $[{1,2,4,3,5,4,2,7,6,5,4,3,1}]$ \\ \hline
$39$ & $w_{13,2}$ & $[{1,4,3,6,5,4,2,7,6,5,4,3,1}]$ \\ \hline
$40$ & $w_{13,3}$ & $[{2,3,1,4,3,5,4,2,6,5,4,3,1}]$ \\ \hline
$41$ & $w_{13,4}$ & $[{2,4,3,6,5,4,2,7,6,5,4,3,1}]$ \\ \hline
$42$ & $w_{13,5}$ & $[{3,1,4,3,5,4,2,7,6,5,4,3,1}]$ \\ \hline
$43$ & $w_{13,6}$ & $[{5,4,3,6,5,4,2,7,6,5,4,3,1}]$ \\ \hline
$44$ & $w_{14,1}$ & $[{1,2,4,3,6,5,4,2,7,6,5,4,3,1}]$ \\ \hline
$45$ & $w_{14,2}$ & $[{1,5,4,3,6,5,4,2,7,6,5,4,3,1}]$ \\ \hline
$46$ & $w_{14,3}$ & $[{2,3,1,4,3,5,4,2,7,6,5,4,3,1}]$ \\ \hline
$47$ & $w_{14,4}$ & $[{2,5,4,3,6,5,4,2,7,6,5,4,3,1}]$ \\ \hline
$48$ & $w_{14,5}$ & $[{3,1,4,3,6,5,4,2,7,6,5,4,3,1}]$ \\ \hline
$49$ & $w_{14,6}$ & $[{4,2,3,1,4,3,5,4,2,6,5,4,3,1}]$ \\ \hline
$50$ & $w_{15,1}$ & $[{1,2,5,4,3,6,5,4,2,7,6,5,4,3,1}]$ \\ \hline
$51$ & $w_{15,2}$ & $[{2,3,1,4,3,6,5,4,2,7,6,5,4,3,1}]$ \\ \hline
$52$ & $w_{15,3}$ & $[{3,1,5,4,3,6,5,4,2,7,6,5,4,3,1}]$ \\ \hline
$53$ & $w_{15,4}$ & $[{4,2,3,1,4,3,5,4,2,7,6,5,4,3,1}]$ \\ \hline
\end{tabular}
}

\vspace{2 mm}

\begin{tabular}{|r|l|r|}
\hline
$54$ & $w_{15,5}$ & $[{4,2,5,4,3,6,5,4,2,7,6,5,4,3,1}]$ \\ \hline
$55$ & $w_{15,6}$ & $[{5,4,2,3,1,4,3,5,4,2,6,5,4,3,1}]$ \\ \hline
$56$ & $w_{16,1}$ & $[{1,4,2,5,4,3,6,5,4,2,7,6,5,4,3,1}]$ \\ \hline
$57$ & $w_{16,2}$ & $[{2,3,1,5,4,3,6,5,4,2,7,6,5,4,3,1}]$ \\ \hline
$58$ & $w_{16,3}$ & $[{3,4,2,5,4,3,6,5,4,2,7,6,5,4,3,1}]$ \\ \hline
$59$ & $w_{16,4}$ & $[{4,2,3,1,4,3,6,5,4,2,7,6,5,4,3,1}]$ \\ \hline
$60$ & $w_{16,5}$ & $[{4,3,1,5,4,3,6,5,4,2,7,6,5,4,3,1}]$ \\ \hline
\end{tabular}

\begin{tabular}{|r|l|r|}
\hline
$61$ & $w_{16,6}$ & $[{5,4,2,3,1,4,3,5,4,2,7,6,5,4,3,1}]$ \\ \hline
$62$ & $w_{16,7}$ & $[{6} ,5,4,2,3,1,4,3,5,4,2,6,5,4,3,1 ]$ \\ \hline
$63$ & $w_{17,1}$ & $[{1,3,4,2,5,4,3,6,5,4,2,7,6,5,4,3,1}]$ \\ \hline
$64$ & $w_{17,2}$ & $[{2,4,3,1,5,4,3,6,5,4,2,7,6,5,4,3,1}]$ \\ \hline
$65$ & $w_{17,3}$ & $[{3,1,4,2,5,4,3,6,5,4,2,7,6,5,4,3,1}]$ \\ \hline
$66$ & $w_{17,4}$ & $[{4,2,3,1,5,4,3,6,5,4,2,7,6,5,4,3,1}]$ \\ \hline
$67$ & $w_{17,5}$ & $[{5,4,2,3,1,4,3,6,5,4,2,7,6,5,4,3,1}]$ \\ \hline
$68$ & $w_{17,6}$ & $[{6,5,4,2,3,1,4,3,5,4,2,7,6,5,4,3,1}]$ \\
\hline
\end{tabular}

{\small \setlength{\tabcolsep}{1 mm} \lineskip 0mm }
\end{center}

{\small 
\begin{tabular}{|r|l|r|}
\hline
$110$ & $w_{25,3}$ & $[{3,1,5,4,2,3,4,5,6,5,4,2,3,1,4,3,5,4,2,7,6,5,4,3,1}] $
\\ \hline
$111$ & $w_{25,4}$ & $[{6,5,4,2,3,4,5,6,7,6} ,5,4,2,3,1,4,3,5,4,2,6,5,4,3,1
] $ \\ \hline
$112$ & $w_{26,1}$ & $[{1,6,5,4,2,3,4,5,6,7,6}
,5,4,2,3,1,4,3,5,4,2,6,5,4,3,1 ]$ \\ \hline
$113$ & $w_{26,2}$ & $[{3,1,5,4,2,3,4,5,6,7,6}
,5,4,2,3,1,4,3,5,4,2,6,5,4,3,1 ]$ \\ \hline
$114$ & $w_{26,3}$ & $[{4,3,1,5,4,2,3,4,5,6,5,4,2,3,1,4,3,5,4,2,7,6,5,4, 3,1}%
]$ \\ \hline
$115$ & $w_{27,1}$ & $[{2,4,3,1,5,4,2,3,4,5,6,5,4,2,3,1,4,3,5,4,2,7,6,5,
4,3,1}]$ \\ \hline
$116$ & $w_{27,2}$ & $[{3,1,6,5,4,2,3,4,5,6,7,6}
,5,4,2,3,1,4,3,5,4,2,6,5,4,3,1 ]$ \\ \hline
$117$ & $w_{27,3}$ & $[{4,3,1,5,4,2,3,4,5,6,7,6}
,5,4,2,3,1,4,3,5,4,2,6,5,4,3,1 ]$ \\ \hline
$118$ & $w_{28,1}$ & $[{2,4,3,1,5,4,2,3,4,5,6,7,6}
,5,4,2,3,1,4,3,5,4,2,6,5,4,3,1 ]$ \\ \hline
$119$ & $w_{28,2}$ & $[{4,3,1,6,5,4,2,3,4,5,6,7,6}
,5,4,2,3,1,4,3,5,4,2,6,5,4,3,1 ]$ \\ \hline
$120$ & $w_{29,1}$ & $[{2,4,3,1,6,5,4,2,3,4,5,6,7,6}
,5,4,2,3,1,4,3,5,4,2,6,5,4,3,1 ]$ \\ \hline
$121$ & $w_{29,2}$ & $[{5,4,3,1,6,5,4,2,3,4,5,6,7,6}
,5,4,2,3,1,4,3,5,4,2,6,5,4,3,1 ]$ \\ \hline
$122$ & $w_{30,1}$ & $[{2,5,4,3,1,6,5,4,2,3,4,5,6,7,6}
,5,4,2,3,1,4,3,5,4,2,6,5,4,3,1 ]$ \\ \hline
$123$ & $w_{31,1}$ & $[{4,2,5,4,3,1,6,5,4,2,3,4,5,6,7,6}
,5,4,2,3,1,4,3,5,4,2,6,5,4,3,1 ]$ \\ \hline
$124$ & $w_{32,1}$ & $[{3,4,2,5,4,3,1,6,5,4,2,3,4,5,6,7,6}
,5,4,2,3,1,4,3,5,4,2,6,5,4,3,1 ]$ \\ \hline
$125$ & $w_{33,1}$ & $[{1,3,4,2,5,4,3,1,6,5,4,2,3,4,5,6,7,6}
,5,4,2,3,1,4,3,5,4,2,6,5,4,3,1 ]$ \\ \hline
\end{tabular}
}

\textbf{6.2. Product with the Euler class }$\omega =s_{1,1}$

\begin{center}
$A_{2}=\left( 
\begin{array}{c}
1%
\end{array}%
\right) $, \ $A_{3}=\left( 
\begin{array}{c}
1%
\end{array}%
\right) $, \ $A_{4}=\left( 
\begin{array}{cc}
1 & 1%
\end{array}%
\right) $, \ $A_{5}=\left( 
\begin{array}{cc}
1 & 0 \\ 
1 & 1%
\end{array}%
\right) $, \ $A_{6}=\left( 
\begin{array}{ccc}
1 & 1 & 0 \\ 
1 & 0 & 1%
\end{array}%
\right) $,

$A_{7}= \left( 
\begin{array}{ccc}
1 & 0 & 1 \\ 
0 & 1 & 1 \\ 
1 & 0 & 0%
\end{array}%
\right) $, \ $A_{8}= \left( 
\begin{array}{cccc}
0 & 0 & 1 & 0 \\ 
1 & 1 & 0 & 0 \\ 
0 & 1 & 1 & 1%
\end{array}%
\right) $, \ $A_{9}= \left( 
\begin{array}{cccc}
1 & 0 & 0 & 0 \\ 
1 & 1 & 1 & 0 \\ 
0 & 1 & 0 & 1 \\ 
0 & 0 & 1 & 1%
\end{array}%
\right) $,

$A_{10}= \left( 
\begin{array}{ccccc}
1 & 1 & 0 & 0 & 0 \\ 
1 & 0 & 1 & 0 & 0 \\ 
0 & 1 & 1 & 1 & 0 \\ 
0 & 0 & 1 & 0 & 1%
\end{array}%
\right) $, \ $A_{11}= \left( 
\begin{array}{ccccc}
1 & 0 & 0 & 0 & 0 \\ 
1 & 1 & 0 & 0 & 0 \\ 
1 & 0 & 0 & 1 & 1 \\ 
0 & 1 & 1 & 0 & 1 \\ 
0 & 0 & 0 & 1 & 0%
\end{array}%
\right) $, \ $A_{12}= \left( 
\begin{array}{cccccc}
0 & 1 & 1 & 0 & 0 & 0 \\ 
1 & 0 & 1 & 0 & 1 & 0 \\ 
1 & 0 & 0 & 1 & 0 & 0 \\ 
0 & 1 & 0 & 0 & 0 & 1 \\ 
0 & 0 & 1 & 1 & 0 & 1%
\end{array}%
\right) $,

$A_{13}= \left( 
\begin{array}{cccccc}
1 & 0 & 1 & 0 & 0 & 0 \\ 
0 & 1 & 0 & 0 & 0 & 0 \\ 
1 & 1 & 0 & 0 & 1 & 0 \\ 
1 & 0 & 0 & 1 & 0 & 0 \\ 
0 & 0 & 1 & 0 & 1 & 0 \\ 
0 & 1 & 0 & 1 & 0 & 1%
\end{array}%
\right) $, \ $A_{14}= \left( 
\begin{array}{cccccc}
1 & 0 & 1 & 0 & 0 & 0 \\ 
1 & 1 & 0 & 0 & 1 & 0 \\ 
0 & 0 & 1 & 0 & 0 & 1 \\ 
1 & 0 & 0 & 1 & 0 & 0 \\ 
0 & 0 & 1 & 0 & 1 & 0 \\ 
0 & 1 & 0 & 1 & 0 & 0%
\end{array}%
\right) $,

$A_{15}= \left( 
\begin{array}{cccccc}
1 & 1 & 0 & 0 & 0 & 0 \\ 
1 & 0 & 1 & 0 & 0 & 0 \\ 
0 & 1 & 0 & 1 & 0 & 0 \\ 
1 & 0 & 0 & 0 & 1 & 0 \\ 
0 & 1 & 1 & 0 & 0 & 0 \\ 
0 & 0 & 0 & 1 & 0 & 1%
\end{array}%
\right) $, \ $A_{16}= \left( 
\begin{array}{ccccccc}
1 & 1 & 0 & 0 & 0 & 0 & 0 \\ 
0 & 1 & 0 & 1 & 0 & 0 & 0 \\ 
0 & 1 & 0 & 0 & 1 & 0 & 0 \\ 
0 & 0 & 0 & 1 & 0 & 1 & 0 \\ 
1 & 0 & 1 & 0 & 0 & 0 & 0 \\ 
0 & 0 & 0 & 0 & 0 & 1 & 1%
\end{array}%
\right) $,

$A_{17}= \left( 
\begin{array}{ccccccc}
1 & 0 & 2 & 1 & 0 & 0 & 0 \\ 
0 & 1 & 1 & 2 & 1 & 0 & 0 \\ 
2 & 0 & 1 & 0 & 0 & 0 & 0 \\ 
0 & 0 & 0 & 1 & 2 & 1 & 0 \\ 
0 & 2 & 0 & 1 & 0 & 0 & 0 \\ 
0 & 0 & 0 & 0 & 1 & 2 & 1 \\ 
0 & 0 & 0 & 0 & 0 & 1 & 2%
\end{array}%
\right) $, \ $A_{18}= \left( 
\begin{array}{cccccc}
1 & 0 & 0 & 0 & 0 & 0 \\ 
0 & 1 & 0 & 0 & 0 & 0 \\ 
1 & 0 & 1 & 0 & 0 & 0 \\ 
0 & 1 & 1 & 1 & 0 & 0 \\ 
0 & 0 & 0 & 1 & 1 & 0 \\ 
0 & 0 & 0 & 0 & 1 & 1 \\ 
0 & 0 & 0 & 0 & 0 & 1%
\end{array}%
\right) $,

$A_{19}= \left( 
\begin{array}{cccccc}
1 & 0 & 0 & 0 & 0 & 0 \\ 
0 & 1 & 1 & 0 & 0 & 0 \\ 
1 & 1 & 0 & 1 & 0 & 0 \\ 
0 & 0 & 1 & 1 & 1 & 0 \\ 
0 & 0 & 0 & 0 & 1 & 1 \\ 
0 & 0 & 0 & 0 & 0 & 1%
\end{array}%
\right) $, \ $A_{20}= \left( 
\begin{array}{cccccc}
1 & 1 & 0 & 0 & 0 & 0 \\ 
1 & 0 & 1 & 0 & 0 & 0 \\ 
0 & 0 & 1 & 1 & 0 & 0 \\ 
0 & 1 & 1 & 0 & 1 & 0 \\ 
0 & 0 & 0 & 1 & 1 & 1 \\ 
0 & 0 & 0 & 0 & 0 & 1%
\end{array}%
\right) $,

{\setlength{\arraycolsep}{1 mm} }
\end{center}

$A_{21}\!\!= \!\!\left( 
\begin{array}{cccccc}
1 & 0 & 0 & 0 & 0 & 0 \\ 
1 & 1 & 0 & 0 & 0 & 0 \\ 
1 & 0 & 1 & 0 & 0 & 1 \\ 
0 & 0 & 1 & 1 & 0 & 0 \\ 
0 & 1 & 1 & 0 & 1 & 0 \\ 
0 & 0 & 0 & 1 & 1 & 0%
\end{array}%
\right) , A_{22}\!\!= \!\!\left( 
\begin{array}{ccccc}
1 & 0 & 1 & 0 & 0 \\ 
1 & 1 & 0 & 0 & 0 \\ 
1 & 0 & 0 & 1 & 1 \\ 
0 & 0 & 0 & 1 & 0 \\ 
0 & 1 & 0 & 1 & 0 \\ 
0 & 0 & 1 & 0 & 1%
\end{array}%
\right) , A_{23}\!\!=\!\! \left( 
\begin{array}{ccccc}
1 & 1 & 0 & 0 & 0 \\ 
1 & 0 & 0 & 0 & 0 \\ 
0 & 1 & 1 & 0 & 0 \\ 
1 & 0 & 0 & 1 & 0 \\ 
0 & 1 & 0 & 1 & 1%
\end{array}%
\right) $,

\begin{center}
$A_{24}= \left( 
\begin{array}{cccc}
1 & 0 & 0 & 0 \\ 
1 & 1 & 1 & 0 \\ 
0 & 0 & 1 & 0 \\ 
1 & 0 & 0 & 1 \\ 
0 & 1 & 0 & 1%
\end{array}%
\right) $, \ $A_{25}= \left( 
\begin{array}{cccc}
1 & 1 & 0 & 0 \\ 
1 & 0 & 1 & 0 \\ 
0 & 1 & 1 & 0 \\ 
1 & 0 & 0 & 1%
\end{array}%
\right) $, \ $A_{26}= \left( 
\begin{array}{ccc}
1 & 1 & 0 \\ 
0 & 1 & 0 \\ 
0 & 1 & 1 \\ 
1 & 0 & 0%
\end{array}%
\right) $,

$A_{27}= \left( 
\begin{array}{ccc}
0 & 1 & 0 \\ 
0 & 1 & 1 \\ 
1 & 0 & 1%
\end{array}%
\right) $, \ $A_{28}= \left( 
\begin{array}{cc}
1 & 0 \\ 
0 & 1 \\ 
1 & 1%
\end{array}%
\right) $, \ $A_{29}= \left( 
\begin{array}{cc}
1 & 0 \\ 
1 & 1%
\end{array}%
\right) $, \ $A_{30}= \left( 
\begin{array}{c}
1 \\ 
1%
\end{array}%
\right) $,

$A_{31}= \left( 
\begin{array}{c}
1%
\end{array}%
\right) $, \ $A_{32}= \left( 
\begin{array}{c}
1%
\end{array}%
\right) $, \ $A_{33}= \left( 
\begin{array}{c}
1%
\end{array}%
\right) $,
\end{center}

\textbf{The additive basis of} $H^{\ast }(E_{7}/D_{6})$ derived from the $%
A_{k}$'s:

\begin{center}
\begin{tabular}{|l|l|l|}
\hline
$H^{0}$ & $Z$ &  \\ \hline
$H^{8}$ & $Z$ & $\overline{s}_{4,1}$ \\ \hline
$H^{12}$ & $Z$ & $\overline{s}_{6,1}$ \\ \hline
$H^{16}$ & $Z$ & $\overline{s}_{8,1}$ \\ \hline
$H^{18}$ & $Z_{2}$ & $\overline{s}_{9,2}$ \\ \hline
$H^{20}$ & $Z$ & $\overline{s}_{10,1}$ \\ \hline
$H^{24}$ & $Z$ & $\overline{s}_{12,2}$ \\ \hline
$H^{26}$ & $Z_{2}$ & $\overline{s}_{13,1}$ \\ \hline
\end{tabular}
\begin{tabular}{|l|l|l|}
\hline
$H^{28}$ & $Z_{3}$ & $\overline{s}_{14,1}$ \\ \hline
$H^{30}$ & $Z_{2}$ & $\overline{s}_{15,1}$ \\ \hline
$H^{32}$ & $Z$ & $\overline{s}_{16,1}$ \\ \hline
$H^{34}$ & $Z_{2}$ & $\overline{s}_{17,2}$ \\ \hline
$H^{38}$ & $Z_{2}$ & $\overline{s}_{19,2}$ \\ \hline
$H^{40}$ & $Z_{3}$ & $\overline{s}_{20,1}$ \\ \hline
$H^{42}$ & $Z_{2}$ & $\overline{s}_{21,3}$ \\ \hline
$H^{50}$ & $Z_{2}$ & $\overline{s}_{25,1}$ \\ \hline
\end{tabular}

\begin{tabular}{|l|l|l|}
\hline
$H^{35}$ & $Z$ & $\beta ^{-1}(s_{17,1}-s_{17,2}-s_{17,3}+s_{17,4}-\
s_{17,5}+s_{17,6}-s_{17,7})$ \\ \hline
$H^{43}$ & $Z$ & $\beta ^{-1}(s_{21,1}-2\,s_{21,2}+s_{21,3}-\
3\,s_{21,4}+2\,s_{21,5}-s_{21,6})$ \\ \hline
$H^{47}$ & $Z$ & $\beta ^{-1}(-2\,s_{23,1}+s_{23,2}-s_{23,3}+\
s_{23,4}-s_{23,5})$ \\ \hline
$H^{51}$ & $Z$ & $\beta ^{-1}(-s_{25,1}+s_{25,2}+s_{25,4})$ \\ \hline
$H^{55}$ & $Z$ & $\beta ^{-1}(-s_{27,1}-s_{27,2}+s_{27,3})$ \\ \hline
$H^{59}$ & $Z$ & $\beta ^{-1}(-s_{29,1}+s_{29,2})$ \\ \hline
$H^{67}$ & $Z$ &  \\ \hline
\end{tabular}
\end{center}

\textbf{6.3. The ring structure on} $H^{even}(E_{7}/D_{6})$:

\begin{tabular}{l}
$\overline{{s}}{_{4,1}}^{2}=\overline{s}_{8,1}$ \\ 
$\overline{s}_{4,1}\,\overline{s}_{6,1}=\overline{s}_{10,1}$ \\ 
$\overline{s}_{4,1}\,\overline{s}_{8,1}=3\,\overline{s}_{12,2}$ \\ 
$\overline{s}_{4,1}\,\overline{s}_{9,2}=3\,\overline{s}_{13,1}$ \\ 
$\overline{s}_{10,1}\,\overline{s}_{4,1}=2\,\overline{s}_{14,1}$ \\ 
$\overline{s}_{12,2}\,\overline{s}_{4,1}=\overline{s}_{16,1}$%
\end{tabular}
\ 
\begin{tabular}{l}
$\overline{s}_{13,1}\,\overline{s}_{4,1}=-3\,\overline{s}_{17,2}$ \\ 
$\overline{s}_{15,1}\,\overline{s}_{4,1}=5\,\overline{s}_{19,2}$ \\ 
$\overline{s}_{16,1}\,\overline{s}_{4,1}=-5\,\overline{s}_{20,1}$ \\ 
$\overline{s}_{17,2}\,\overline{s}_{4,1}=-\overline{s}_{21,3}$ \\ 
$\overline{s}_{21,3}\,\overline{s}_{4,1}=5\,\overline{s}_{25,1}$ \\ 
$\overline{{s}}{_{6,1}}^{2}=\overline{s}_{12,2}$%
\end{tabular}
\ 
\begin{tabular}{l}
$\overline{s}_{6,1}\,\overline{s}_{8,1}=2\,\overline{s}_{14,1}$ \\ 
$\overline{s}_{6,1}\,\overline{s}_{9,2}=5\,\overline{s}_{15,1}$ \\ 
$\overline{s}_{10,1}\,\overline{s}_{6,1}=\overline{s}_{16,1}$ \\ 
$\overline{s}_{13,1}\,\overline{s}_{6,1}=5\,\overline{s}_{19,2}$ \\ 
$\overline{s}_{14,1}\,\overline{s}_{6,1}=-10\,\overline{s}_{20,1}$%
\end{tabular}

\bigskip

\textbf{6.4. The structure matrices }$M(\pi _{m})$\textbf{\ }are with
respect to the ordered monomial base $B(2m)$ $:$

$%
%
%
\right) $ }

\textbf{6.5. The nullspaces of }$M(\pi _{m})$ for $m=9,12,14,18$ are
respectively

\begin{center}
$N(\pi _{9})=\left( 
\begin{array}{ccccc}
-2 & -3 & -4 & -2 & 2%
\end{array}%
\right) $;$N(\pi _{12})=\left( 
\begin{array}{cccccccc}
3 & -1 & 0 & -2 & -6 & -6 & 0 & 2 \\ 
3 & -1 & 0 & 0 & -3 & -2 & 2 & 0%
\end{array}%
\right) $;

$N(\pi _{14})=\left( 
\begin{array}{cccccccccc}
-3 & 0 & -3 & -6 & -6 & -2 & 0 & 0 & 0 & 1 \\ 
-12 & 6 & 0 & -19 & -12 & -4 & 0 & 0 & 6 & 0 \\ 
6 & -12 & -9 & -1 & -12 & 2 & 0 & 6 & 0 & 0 \\ 
-6 & 6 & 0 & -5 & 0 & -2 & 3 & 0 & 0 & 0%
\end{array}%
\right) $;

{\small \setlength{\arraycolsep}{0.6 mm} $N(\pi _{18})=\left( 
\begin{array}{ccccccccccccccccc}
-608 & 836 & -608 & -1140 & -1140 & -361 & 0 & 0 & 0 & 0 & 0 & 0 & 0 & 0 & \
0 & 0 & 28 \\ 
-110 & 174 & -166 & -138 & -180 & -85 & 0 & 0 & 0 & 0 & 0 & 0 & 0 & 0 & 0 & 
14 & 0 \\ 
-236 & 244 & -180 & -432 & -600 & -29 & 0 & 0 & 0 & 0 & 0 & 0 & 0 & 0 & 56\ 
& 0 & 0 \\ 
-55 & 87 & -104 & -48 & -69 & -67 & 0 & 0 & 0 & 0 & 0 & 0 & 0 & 21 & 0 & 0\ 
& 0 \\ 
-3 & 5 & 8 & -24 & -3 & -2 & 0 & 0 & 0 & 0 & 0 & 0 & 2 & 0 & 0 & 0 & 0 \\ 
-166 & 258 & -236 & -138 & -348 & -1 & 0 & 0 & 0 & 0 & 0 & 84 & 0 & 0 & 0 & 
\ 0 & 0 \\ 
-6 & 10 & -20 & -6 & -6 & -25 & 0 & 0 & 0 & 0 & 14 & 0 & 0 & 0 & 0 & 0 & \ 0
\\ 
-4 & -12 & 4 & -12 & -24 & -1 & 0 & 0 & 0 & 12 & 0 & 0 & 0 & 0 & 0 & 0 & \ 0
\\ 
-18 & 30 & 24 & -102 & -18 & -5 & 0 & 0 & 28 & 0 & 0 & 0 & 0 & 0 & 0 & 0 & \
0 \\ 
-23 & 36 & -37 & -9 & -30 & 1 & 0 & 21 & 0 & 0 & 0 & 0 & 0 & 0 & 0 & 0 & \ 0
\\ 
4 & 40 & 88 & -192 & -24 & -23 & 56 & 0 & 0 & 0 & 0 & 0 & 0 & 0 & 0 & 0 & \ 0%
\end{array}%
\right) $ }
\end{center}

\section{$(G,H_{{s}})=(E_{8},E_{7})$}

\textbf{7.1. The minimized decompositions of elements in }$W(E_{7}\cdot
S^{1},E_{8})$:

%% all cells
\begin{tabular}{|l|l|l|}
\hline
& $w_{i,j}$ & decomposition \\ \hline
$1$ & $w_{1,1}$ & $[{8}]$ \\ 
$2$ & $w_{2,1}$ & $[{7, 8}]$ \\ 
$3$ & $w_{3,1}$ & $[{6, 7, 8}]$ \\ 
$4$ & $w_{4,1}$ & $[{5, 6, 7, 8}]$ \\ 
$5$ & $w_{5,1}$ & $[{4, 5, 6, 7, 8}]$ \\ 
$6$ & $w_{6,1}$ & $[{2, 4, 5, 6, 7, 8}]$ \\ 
$7$ & $w_{6,2}$ & $[{3, 4, 5, 6, 7, 8}]$ \\ 
$8$ & $w_{7,1}$ & $[{1, 3, 4, 5, 6, 7, 8}]$ \\ 
$9$ & $w_{7,2}$ & $[{2, 3, 4, 5, 6, 7, 8}]$ \\ 
$10$ & $w_{8,1}$ & $[{1, 2, 3, 4, 5, 6, 7, 8}]$ \\ 
$11$ & $w_{8,2}$ & $[{4, 2, 3, 4, 5, 6, 7, 8}]$ \\ 
$12$ & $w_{9,1}$ & $[{1, 4, 2, 3, 4, 5, 6, 7, 8}]$ \\ 
$13$ & $w_{9,2}$ & $[{5, 4, 2, 3, 4, 5, 6, 7, 8}]$ \\ 
$14$ & $w_{10,1}$ & $[{1, 5, 4, 2, 3, 4, 5, 6, 7, 8}]$ \\ 
$15$ & $w_{10,2}$ & $[{3, 1, 4, 2, 3, 4, 5, 6, 7, 8}]$ \\ 
$16$ & $w_{10,3}$ & $[{6, 5, 4, 2, 3, 4, 5, 6, 7, 8}]$ \\ 
$17$ & $w_{11,1}$ & $[{1, 6, 5, 4, 2, 3, 4, 5, 6, 7, 8}]$ \\ 
$18$ & $w_{11,2}$ & $[{3, 1, 5, 4, 2, 3, 4, 5, 6, 7, 8}]$ \\ 
$19$ & $w_{11,3}$ & $[{7, 6, 5, 4, 2, 3, 4, 5, 6, 7, 8}]$ \\ 
$20$ & $w_{12,1}$ & $[{1, 7, 6, 5, 4, 2, 3, 4, 5, 6, 7, 8}]$ \\ 
$21$ & $w_{12,2}$ & $[{3, 1, 6, 5, 4, 2, 3, 4, 5, 6, 7, 8}]$ \\ 
$22$ & $w_{12,3}$ & $[{4, 3, 1, 5, 4, 2, 3, 4, 5, 6, 7, 8}]$ \\ 
$23$ & $w_{12,4}$ & $[{8, 7, 6, 5, 4, 2, 3, 4, 5, 6, 7, 8}]$ \\ 
$24$ & $w_{13,1}$ & $[{1, 8, 7, 6, 5, 4, 2, 3, 4, 5, 6, 7, 8}]$ \\ 
$25$ & $w_{13,2}$ & $[{2, 4, 3, 1, 5, 4, 2, 3, 4, 5, 6, 7, 8}]$ \\ 
$26$ & $w_{13,3}$ & $[{3, 1, 7, 6, 5, 4, 2, 3, 4, 5, 6, 7, 8}]$ \\ 
$27$ & $w_{13,4}$ & $[{4, 3, 1, 6, 5, 4, 2, 3, 4, 5, 6, 7, 8}]$ \\ 
$28$ & $w_{14,1}$ & $[{2, 4, 3, 1, 6, 5, 4, 2, 3, 4, 5, 6, 7, 8}]$ \\ 
$29$ & $w_{14,2}$ & $[{3, 1, 8, 7, 6, 5, 4, 2, 3, 4, 5, 6, 7, 8}]$ \\ 
$30$ & $w_{14,3}$ & $[{4, 3, 1, 7, 6, 5, 4, 2, 3, 4, 5, 6, 7, 8}]$ \\ 
$31$ & $w_{14,4}$ & $[{5, 4, 3, 1, 6, 5, 4, 2, 3, 4, 5, 6, 7, 8}]$ \\ 
$32$ & $w_{15,1}$ & $[{2, 4, 3, 1, 7, 6, 5, 4, 2, 3, 4, 5, 6, 7, 8}]$ \\ 
$33$ & $w_{15,2}$ & $[{2, 5, 4, 3, 1, 6, 5, 4, 2, 3, 4, 5, 6, 7, 8}]$ \\ 
\hline
\end{tabular}

\begin{tabular}{|l|l|l|}
\hline
$34$ & $w_{15,3}$ & $[{4, 3, 1, 8, 7, 6, 5, 4, 2, 3, 4, 5, 6, 7, 8}]$ \\ 
$35$ & $w_{15,4}$ & $[{5, 4, 3, 1, 7, 6, 5, 4, 2, 3, 4, 5, 6, 7, 8}]$ \\ 
$36$ & $w_{16,1}$ & $[{2, 4, 3, 1, 8, 7, 6, 5, 4, 2, 3, 4, 5, 6, 7, 8}]$ \\ 
$37$ & $w_{16,2}$ & $[{2, 5, 4, 3, 1, 7, 6, 5, 4, 2, 3, 4, 5, 6, 7, 8}]$ \\ 
$38$ & $w_{16,3}$ & $[{4, 2, 5, 4, 3, 1, 6, 5, 4, 2, 3, 4, 5, 6, 7, 8}]$ \\ 
$39$ & $w_{16,4}$ & $[{5, 4, 3, 1, 8, 7, 6, 5, 4, 2, 3, 4, 5, 6, 7, 8}]$ \\ 
$40$ & $w_{16,5}$ & $[{6, 5, 4, 3, 1, 7, 6, 5, 4, 2, 3, 4, 5, 6, 7, 8}]$ \\ 
$41$ & $w_{17,1}$ & $[{2, 5, 4, 3, 1, 8, 7, 6, 5, 4, 2, 3, 4, 5, 6, 7, 8}]$
\\
\hline
\end{tabular}

\begin{tabular}{|l|l|l|}
\hline
$84$ & $w_{24,1}$ & $[{1, 2, 4, 3, 5, 4, 2, 6, 5, 4, 3, 1, 8, 7, 6, 5, 4, 2,
3, 4, 5, 6, 7, \ 8}]$ \\ 
$85$ & $w_{24,2}$ & $[{1, 3, 6, 5, 4, 2, 7, 6, 5, 4, 3, 1, 8, 7, 6, 5, 4, 2,
3, 4, 5, 6, 7, \ 8}]$ \\ 
$86$ & $w_{24,3}$ & $[{1, 4, 3, 5, 4, 2, 7, 6, 5, 4, 3, 1, 8, 7, 6, 5, 4, 2,
3, 4, 5, 6, 7, \ 8}]$ \\ 
$87$ & $w_{24,4}$ & $[{2, 3, 1, 4, 3, 5, 4, 2, 6, 5, 4, 3, 1, 7, 6, 5, 4, 2,
3, 4, 5, 6, 7, \ 8}]$ \\ 
$88$ & $w_{24,5}$ & $[{2, 4, 3, 5, 4, 2, 7, 6, 5, 4, 3, 1, 8, 7, 6, 5, 4, 2,
3, 4, 5, 6, 7, \ 8}]$ \\ 
$89$ & $w_{24,6}$ & $[{3, 1, 4, 3, 5, 4, 2, 6, 5, 4, 3, 1, 8, 7, 6, 5, 4, 2,
3, 4, 5, 6, 7, \ 8}]$ \\ 
$90$ & $w_{24,7}$ & $[{4, 3, 6, 5, 4, 2, 7, 6, 5, 4, 3, 1, 8, 7, 6, 5, 4, 2,
3, 4, 5, 6, 7, \ 8}]$ \\ 
$91$ & $w_{25,1}$ & $[{1, 2, 4, 3, 5, 4, 2, 7, 6, 5, 4, 3, 1, 8, 7, 6, 5, 4,
2, 3, 4, 5, 6, 7, \ 8}]$ \\
\hline
\end{tabular}

\begin{tabular}{|l|l|l|}
\hline
$135$ & $w_{31,1}$ & $[{1, 3, 4, 2, 3, 1, 5, 4, 3, 6, 5, 4, 2, 7, 6, 5, 4,
3, 1, 8, 7, 6, 5, 4, \ 2, 3, 4, 5, 6, 7, 8}]$ \\ 
$136$ & $w_{31,2}$ & $[{2, 3, 4, 2, 3, 1, 5, 4, 3, 6, 5, 4, 2, 7, 6, 5, 4,
3, 1, 8, 7, 6, 5, 4, \ 2, 3, 4, 5, 6, 7, 8}]$ \\ 
$137$ & $w_{31,3}$ & $[{2, 4, 5, 4, 2, 3, 1, 4, 3, 6, 5, 4, 2, 7, 6, 5, 4,
3, 1, 8, 7, 6, 5, 4, \ 2, 3, 4, 5, 6, 7, 8}]$ \\ 
$138$ & $w_{31,4}$ & $[{3, 4, 5, 4, 2, 3, 1, 4, 3, 6, 5, 4, 2, 7, 6, 5, 4,
3, 1, 8, 7, 6, 5, 4, \ 2, 3, 4, 5, 6, 7, 8}]$ \\ 
$139$ & $w_{31,5}$ & $[{4, 5, 6, 5, 4, 2, 3, 1, 4, 3, 5, 4, 2, 7, 6, 5, 4,
3, 1, 8, 7, 6, 5, 4, \ 2, 3, 4, 5, 6, 7, 8}]$ \\ 
$140$ & $w_{31,6}$ & $[{5, 6, 7, 6, 5, 4, 2, 3, 1, 4, 3, 5, 4, 2, 6, 5, 4,
3, 1, 8, 7, 6, 5, 4, \ 2, 3, 4, 5, 6, 7, 8}]$ \\ 
$141$ & $w_{31,7}$ & $[{6, 7, 8, 7, 6, 5, 4, 2, 3, 1, 4, 3, 5, 4, 2, 6, 5,
4, 3, 1, 7, 6, 5, 4, \ 2, 3, 4, 5, 6, 7, 8}]$ \\ 
$142$ & $w_{32,1}$ & $[{1, 2, 3, 4, 2, 3, 1, 5, 4, 3, 6, 5, 4, 2, 7, 6, 5,
4, 3, 1, 8, 7, 6, 5, \ 4, 2, 3, 4, 5, 6, 7, 8}]$ \\
\hline
\end{tabular}

\begin{tabular}{|l|l|l|}
\hline
$185$ & $w_{38,4}$ & $[{3, 1, 5, 4, 2, 3, 4, 5, 6, 7, 6, 5, 4, 2, 3, 1, 4,
3, 5, 4, 2, 6, 5, 4, \ 3, 1, 8, 7, 6, 5, 4, 2, 3, 4, 5, 6, 7, 8}]$ \\ 
$186$ & $w_{38,5}$ & $[{4, 3, 1, 5, 4, 2, 3, 4, 5, 6, 5, 4, 2, 3, 1, 4, 3,
5, 4, 2, 7, 6, 5, 4, \ 3, 1, 8, 7, 6, 5, 4, 2, 3, 4, 5, 6, 7, 8}]$ \\ 
$187$ & $w_{38,6}$ & $[{6, 5, 4, 2, 3, 4, 5, 6, 7, 8, 7, 6, 5, 4, 2, 3, 1,
4, 3, 5, 4, 2, 6, 5, \ 4, 3, 1, 7, 6, 5, 4, 2, 3, 4, 5, 6, 7, 8}]$ \\ 
$188$ & $w_{39,1}$ & $[{1, 6, 5, 4, 2, 3, 4, 5, 6, 7, 8, 7, 6, 5, 4, 2, 3,
1, 4, 3, 5, 4, 2, 6, \ 5, 4, 3, 1, 7, 6, 5, 4, 2, 3, 4, 5, 6, 7, 8}]$ \\ 
$189$ & $w_{39,2}$ & $[{2, 4, 3, 1, 5, 4, 2, 3, 4, 5, 6, 5, 4, 2, 3, 1, 4,
3, 5, 4, 2, 7, 6, 5, \ 4, 3, 1, 8, 7, 6, 5, 4, 2, 3, 4, 5, 6, 7, 8}]$ \\ 
$190$ & $w_{39,3}$ & $[{3, 1, 5, 4, 2, 3, 4, 5, 6, 7, 8, 7, 6, 5, 4, 2, 3,
1, 4, 3, 5, 4, 2, 6, \ 5, 4, 3, 1, 7, 6, 5, 4, 2, 3, 4, 5, 6, 7, 8}]$ \\ 
$191$ & $w_{39,4}$ & $[{3, 1, 6, 5, 4, 2, 3, 4, 5, 6, 7, 6, 5, 4, 2, 3, 1,
4, 3, 5, 4, 2, 6, 5, \ 4, 3, 1, 8, 7, 6, 5, 4, 2, 3, 4, 5, 6, 7, 8}]$ \\ 
$192$ & $w_{39,5}$ & $[{4, 3, 1, 5, 4, 2, 3, 4, 5, 6, 7, 6, 5, 4, 2, 3, 1,
4, 3, 5, 4, 2, 6, 5, \ 4, 3, 1, 8, 7, 6, 5, 4, 2, 3, 4, 5, 6, 7, 8}]$ \\
\hline
\end{tabular}

\begin{tabular}{|l|l|l|}
\hline
$236$ & $w_{54,1}$ & $[{5,4,2,3,1,4,3,5,4,2,6,5,4,3,1,7,6,5,4,2,3,4,5,6,\
7,8,7,6,5,4,2,3,1,4,3,5,4,2,6,5,4,3,1,7,6,5,4,2,3,\ 4,5,6,7,8}]$ \\ 
$237$ & $w_{55,1}$ & $[{6,5,4,2,3,1,4,3,5,4,2,6,5,4,3,1,7,6,5,4,2,3,4,5,\
6,7,8,7,6,5,4,2,3,1,4,3,5,4,2,6,5,4,3,1,7,6,5,4,2,\ 3,4,5,6,7,8}]$ \\ 
$238$ & $w_{56,1}$ & $[{7,6,5,4,2,3,1,4,3,5,4,2,6,5,4,3,1,7,6,5,4,2,3,4,\
5,6,7,8,7,6,5,4,2,3,1,4,3,5,4,2,6,5,4,3,1,7,6,5,4,\ 2,3,4,5,6,7,8}]$ \\ 
$239$ & $w_{57,1}$ & $[{8,7,6,5,4,2,3,1,4,3,5,4,2,6,5,4,3,1,7,6,5,4,2,3,\
4,5,6,7,8,7,6,5,4,2,3,1,4,3,5,4,2,6,5,4,3,1,7,6,5,\ 4,2,3,4,5,6,7,8}]$ \\ 
\hline
\end{tabular}

\bigskip

\textbf{7.2. Product with the Euler class }$\omega =s_{{1,1}}$

$A_{2} = \left( 
\begin{array}{c}
1%
\end{array}%
\right) $ $A_{3} = \left( 
\begin{array}{c}
1%
\end{array}%
\right) $ $A_{4} = \left( 
\begin{array}{c}
1%
\end{array}%
\right) $ $A_{5} = \left( 
\begin{array}{c}
1%
\end{array}%
\right) $ $A_{6} = \left( 
\begin{array}{cc}
1 & 1%
\end{array}%
\right) $

$A_{7} = \left( 
\begin{array}{cc}
0 & 1 \\ 
1 & 1%
\end{array}%
\right) $ $A_{8} = \left( 
\begin{array}{cc}
1 & 0 \\ 
1 & 1%
\end{array}%
\right) $ $A_{9} = \left( 
\begin{array}{cc}
1 & 0 \\ 
1 & 1%
\end{array}%
\right) $ $A_{10} = \left( 
\begin{array}{ccc}
1 & 1 & 0 \\ 
1 & 0 & 1%
\end{array}%
\right) $

$A_{11} = \left( 
\begin{array}{ccc}
1 & 1 & 0 \\ 
0 & 1 & 0 \\ 
1 & 0 & 1%
\end{array}%
\right) $ $A_{12} = \left( 
\begin{array}{cccc}
1 & 1 & 0 & 0 \\ 
0 & 1 & 1 & 0 \\ 
1 & 0 & 0 & 1%
\end{array}%
\right) $ $A_{13} = \left( 
\begin{array}{cccc}
1 & 0 & 1 & 0 \\ 
0 & 0 & 1 & 1 \\ 
0 & 1 & 0 & 1 \\ 
1 & 0 & 0 & 0%
\end{array}%
\right) $ $A_{14} = \left( 
\begin{array}{cccc}
0 & 1 & 0 & 0 \\ 
1 & 0 & 0 & 0 \\ 
0 & 1 & 1 & 0 \\ 
1 & 0 & 1 & 1%
\end{array}%
\right) $ $A_{15} = \left( 
\begin{array}{cccc}
1 & 1 & 0 & 0 \\ 
0 & 0 & 1 & 0 \\ 
1 & 0 & 1 & 1 \\ 
0 & 1 & 0 & 1%
\end{array}%
\right) $ $A_{16} = \left( 
\begin{array}{ccccc}
1 & 1 & 0 & 0 & 0 \\ 
0 & 1 & 1 & 0 & 0 \\ 
1 & 0 & 0 & 1 & 0 \\ 
0 & 1 & 0 & 1 & 1%
\end{array}%
\right) $ $A_{17} = \left( 
\begin{array}{ccccc}
1 & 0 & 0 & 0 & 0 \\ 
1 & 1 & 0 & 1 & 0 \\ 
0 & 0 & 1 & 1 & 0 \\ 
1 & 0 & 0 & 0 & 1 \\ 
0 & 1 & 0 & 0 & 1%
\end{array}%
\right) $ $A_{18} = \left( 
\begin{array}{cccccc}
0 & 1 & 0 & 1 & 0 & 0 \\ 
0 & 1 & 0 & 0 & 1 & 0 \\ 
1 & 0 & 1 & 0 & 0 & 0 \\ 
0 & 0 & 1 & 1 & 1 & 0 \\ 
0 & 1 & 0 & 0 & 0 & 1%
\end{array}%
\right) $ $A_{19} = \left( 
\begin{array}{cccccc}
1 & 0 & 0 & 0 & 0 & 0 \\ 
0 & 1 & 0 & 0 & 1 & 0 \\ 
1 & 0 & 1 & 1 & 0 & 0 \\ 
0 & 0 & 1 & 0 & 1 & 0 \\ 
0 & 0 & 0 & 1 & 1 & 1 \\ 
0 & 1 & 0 & 0 & 0 & 0%
\end{array}%
\right) $

$A_{20} = \left( 
\begin{array}{cccccc}
1 & 1 & 0 & 0 & 0 & 0 \\ 
0 & 0 & 0 & 0 & 1 & 0 \\ 
1 & 0 & 1 & 0 & 0 & 0 \\ 
0 & 1 & 1 & 1 & 0 & 0 \\ 
0 & 0 & 1 & 0 & 1 & 1 \\ 
0 & 0 & 0 & 1 & 0 & 1%
\end{array}%
\right) $ $A_{21} = \left( 
\begin{array}{cccccc}
1 & 0 & 0 & 0 & 0 & 0 \\ 
1 & 1 & 0 & 0 & 0 & 0 \\ 
1 & 0 & 1 & 1 & 0 & 0 \\ 
0 & 1 & 0 & 1 & 1 & 0 \\ 
0 & 0 & 1 & 0 & 0 & 1 \\ 
0 & 0 & 0 & 1 & 0 & 1%
\end{array}%
\right) $

$A_{22} = \left( 
\begin{array}{ccccccc}
1 & 1 & 0 & 0 & 0 & 0 & 0 \\ 
0 & 1 & 1 & 0 & 0 & 0 & 0 \\ 
1 & 0 & 0 & 0 & 1 & 0 & 0 \\ 
0 & 1 & 0 & 0 & 1 & 1 & 0 \\ 
0 & 0 & 1 & 1 & 0 & 1 & 0 \\ 
0 & 0 & 0 & 0 & 1 & 0 & 1%
\end{array}%
\right) $ $A_{23} = \left( 
\begin{array}{ccccccc}
0 & 1 & 0 & 0 & 0 & 0 & 0 \\ 
0 & 1 & 1 & 0 & 0 & 0 & 0 \\ 
1 & 0 & 1 & 0 & 1 & 0 & 0 \\ 
1 & 0 & 0 & 1 & 0 & 0 & 0 \\ 
0 & 1 & 0 & 0 & 0 & 1 & 1 \\ 
0 & 0 & 1 & 1 & 0 & 0 & 1 \\ 
0 & 0 & 0 & 0 & 0 & 1 & 0%
\end{array}%
\right) $

$A_{24} = \left( 
\begin{array}{ccccccc}
1 & 0 & 0 & 1 & 0 & 0 & 0 \\ 
0 & 1 & 1 & 0 & 0 & 0 & 0 \\ 
1 & 0 & 1 & 0 & 0 & 1 & 0 \\ 
1 & 0 & 0 & 0 & 1 & 0 & 0 \\ 
0 & 0 & 0 & 1 & 0 & 1 & 0 \\ 
0 & 1 & 0 & 0 & 0 & 0 & 1 \\ 
0 & 0 & 1 & 0 & 1 & 0 & 1%
\end{array}%
\right) $ $A_{25} = \left( 
\begin{array}{ccccccc}
1 & 0 & 1 & 0 & 0 & 0 & 0 \\ 
0 & 1 & 0 & 0 & 0 & 0 & 0 \\ 
1 & 1 & 0 & 0 & 1 & 0 & 0 \\ 
0 & 0 & 1 & 0 & 0 & 1 & 0 \\ 
1 & 0 & 0 & 1 & 0 & 0 & 0 \\ 
0 & 0 & 1 & 0 & 1 & 0 & 0 \\ 
0 & 1 & 0 & 1 & 0 & 0 & 1%
\end{array}%
\right) $

$A_{26} = \left( 
\begin{array}{ccccccc}
1 & 0 & 1 & 0 & 0 & 0 & 0 \\ 
1 & 1 & 0 & 0 & 1 & 0 & 0 \\ 
0 & 0 & 1 & 0 & 0 & 1 & 0 \\ 
1 & 0 & 0 & 1 & 0 & 0 & 0 \\ 
0 & 0 & 1 & 0 & 1 & 0 & 0 \\ 
0 & 0 & 0 & 0 & 0 & 1 & 1 \\ 
0 & 1 & 0 & 1 & 0 & 0 & 0%
\end{array}%
\right) $ $A_{27} = \left( 
\begin{array}{ccccccc}
1 & 1 & 0 & 0 & 0 & 0 & 0 \\ 
1 & 0 & 1 & 0 & 0 & 0 & 0 \\ 
0 & 1 & 0 & 1 & 0 & 0 & 0 \\ 
1 & 0 & 0 & 0 & 1 & 0 & 0 \\ 
0 & 1 & 1 & 0 & 0 & 0 & 0 \\ 
0 & 0 & 0 & 1 & 0 & 1 & 0 \\ 
0 & 0 & 0 & 0 & 0 & 1 & 1%
\end{array}%
\right) $

$A_{28} = \left( 
\begin{array}{cccccccc}
1 & 1 & 0 & 0 & 0 & 0 & 0 & 0 \\ 
0 & 1 & 0 & 1 & 0 & 0 & 0 & 0 \\ 
0 & 1 & 0 & 0 & 1 & 0 & 0 & 0 \\ 
0 & 0 & 0 & 1 & 0 & 1 & 0 & 0 \\ 
1 & 0 & 1 & 0 & 0 & 0 & 0 & 0 \\ 
0 & 0 & 0 & 0 & 0 & 1 & 1 & 0 \\ 
0 & 0 & 0 & 0 & 0 & 0 & 1 & 1%
\end{array}%
\right) $ $A_{29} = \left( 
\begin{array}{cccccccc}
1 & 0 & 2 & 1 & 0 & 0 & 0 & 0 \\ 
0 & 1 & 1 & 2 & 1 & 0 & 0 & 0 \\ 
2 & 0 & 1 & 0 & 0 & 0 & 0 & 0 \\ 
0 & 0 & 0 & 1 & 2 & 1 & 0 & 0 \\ 
0 & 2 & 0 & 1 & 0 & 0 & 0 & 0 \\ 
0 & 0 & 0 & 0 & 1 & 2 & 1 & 0 \\ 
0 & 0 & 0 & 0 & 0 & 1 & 2 & 1 \\ 
0 & 0 & 0 & 0 & 0 & 0 & 1 & 2%
\end{array}%
\right) $

$A_{30} = \left( 
\begin{array}{ccccccc}
1 & 0 & 0 & 0 & 0 & 0 & 0 \\ 
0 & 1 & 0 & 0 & 0 & 0 & 0 \\ 
1 & 0 & 1 & 0 & 0 & 0 & 0 \\ 
0 & 1 & 1 & 1 & 0 & 0 & 0 \\ 
0 & 0 & 0 & 1 & 1 & 0 & 0 \\ 
0 & 0 & 0 & 0 & 1 & 1 & 0 \\ 
0 & 0 & 0 & 0 & 0 & 1 & 1 \\ 
0 & 0 & 0 & 0 & 0 & 0 & 1%
\end{array}%
\right) $ $A_{31} = \left( 
\begin{array}{ccccccc}
1 & 0 & 0 & 0 & 0 & 0 & 0 \\ 
0 & 1 & 1 & 0 & 0 & 0 & 0 \\ 
1 & 1 & 0 & 1 & 0 & 0 & 0 \\ 
0 & 0 & 1 & 1 & 1 & 0 & 0 \\ 
0 & 0 & 0 & 0 & 1 & 1 & 0 \\ 
0 & 0 & 0 & 0 & 0 & 1 & 1 \\ 
0 & 0 & 0 & 0 & 0 & 0 & 1%
\end{array}%
\right) $

$A_{32} = \left( 
\begin{array}{ccccccc}
1 & 1 & 0 & 0 & 0 & 0 & 0 \\ 
1 & 0 & 1 & 0 & 0 & 0 & 0 \\ 
0 & 0 & 1 & 1 & 0 & 0 & 0 \\ 
0 & 1 & 1 & 0 & 1 & 0 & 0 \\ 
0 & 0 & 0 & 1 & 1 & 1 & 0 \\ 
0 & 0 & 0 & 0 & 0 & 1 & 1 \\ 
0 & 0 & 0 & 0 & 0 & 0 & 1%
\end{array}%
\right) $ $A_{33} = \left( 
\begin{array}{ccccccc}
1 & 0 & 0 & 0 & 0 & 0 & 0 \\ 
1 & 1 & 0 & 0 & 0 & 0 & 0 \\ 
1 & 0 & 1 & 0 & 0 & 1 & 0 \\ 
0 & 0 & 1 & 1 & 0 & 0 & 0 \\ 
0 & 1 & 1 & 0 & 1 & 0 & 0 \\ 
0 & 0 & 0 & 1 & 1 & 0 & 1 \\ 
0 & 0 & 0 & 0 & 0 & 0 & 1%
\end{array}%
\right) $

$A_{34} = \left( 
\begin{array}{ccccccc}
1 & 0 & 1 & 0 & 0 & 0 & 0 \\ 
1 & 1 & 0 & 0 & 0 & 0 & 0 \\ 
1 & 0 & 0 & 1 & 0 & 0 & 1 \\ 
0 & 0 & 0 & 1 & 1 & 0 & 0 \\ 
0 & 1 & 0 & 1 & 0 & 1 & 0 \\ 
0 & 0 & 1 & 0 & 0 & 0 & 1 \\ 
0 & 0 & 0 & 0 & 1 & 1 & 0%
\end{array}%
\right) $ $A_{35} = \left( 
\begin{array}{ccccccc}
1 & 0 & 1 & 0 & 0 & 0 & 0 \\ 
1 & 1 & 0 & 0 & 0 & 0 & 0 \\ 
0 & 0 & 1 & 0 & 1 & 0 & 0 \\ 
1 & 0 & 0 & 1 & 0 & 1 & 0 \\ 
0 & 0 & 0 & 1 & 0 & 0 & 0 \\ 
0 & 1 & 0 & 1 & 0 & 0 & 0 \\ 
0 & 0 & 1 & 0 & 0 & 1 & 1%
\end{array}%
\right) $

$A_{36} = \left( 
\begin{array}{cccccc}
1 & 1 & 0 & 0 & 0 & 0 \\ 
1 & 0 & 0 & 0 & 0 & 0 \\ 
0 & 1 & 1 & 1 & 0 & 0 \\ 
1 & 0 & 0 & 0 & 1 & 0 \\ 
0 & 0 & 0 & 1 & 0 & 0 \\ 
0 & 1 & 0 & 0 & 1 & 1 \\ 
0 & 0 & 1 & 0 & 0 & 1%
\end{array}%
\right) $ $A_{37} = \left( 
\begin{array}{cccccc}
1 & 0 & 0 & 0 & 0 & 0 \\ 
1 & 1 & 1 & 0 & 0 & 0 \\ 
0 & 1 & 0 & 1 & 0 & 0 \\ 
0 & 0 & 1 & 1 & 0 & 0 \\ 
1 & 0 & 0 & 0 & 1 & 0 \\ 
0 & 1 & 0 & 0 & 1 & 1%
\end{array}%
\right) $

$A_{38} = \left( 
\begin{array}{cccccc}
1 & 0 & 1 & 0 & 0 & 0 \\ 
1 & 1 & 0 & 1 & 0 & 0 \\ 
0 & 0 & 1 & 1 & 0 & 0 \\ 
0 & 0 & 0 & 1 & 1 & 0 \\ 
1 & 0 & 0 & 0 & 0 & 1 \\ 
0 & 1 & 0 & 0 & 0 & 1%
\end{array}%
\right) $ $A_{39} = \left( 
\begin{array}{cccccc}
1 & 0 & 1 & 0 & 0 & 0 \\ 
1 & 0 & 0 & 1 & 0 & 0 \\ 
0 & 0 & 1 & 0 & 0 & 0 \\ 
0 & 0 & 1 & 1 & 1 & 0 \\ 
0 & 1 & 0 & 0 & 1 & 0 \\ 
1 & 0 & 0 & 0 & 0 & 1%
\end{array}%
\right) $

$A_{40} = \left( 
\begin{array}{ccccc}
1 & 0 & 1 & 0 & 0 \\ 
0 & 1 & 0 & 0 & 0 \\ 
0 & 0 & 1 & 1 & 0 \\ 
0 & 0 & 1 & 0 & 1 \\ 
0 & 1 & 0 & 1 & 1 \\ 
1 & 0 & 0 & 0 & 0%
\end{array}%
\right) $ $A_{41} = \left( 
\begin{array}{ccccc}
0 & 0 & 1 & 0 & 0 \\ 
1 & 1 & 0 & 0 & 0 \\ 
0 & 0 & 1 & 1 & 0 \\ 
1 & 0 & 0 & 1 & 0 \\ 
0 & 1 & 0 & 1 & 1%
\end{array}%
\right) $ $A_{42} = \left( 
\begin{array}{cccc}
1 & 0 & 0 & 0 \\ 
1 & 1 & 0 & 0 \\ 
0 & 0 & 1 & 0 \\ 
1 & 0 & 1 & 1 \\ 
0 & 1 & 0 & 1%
\end{array}%
\right) $

$A_{43} = \left( 
\begin{array}{cccc}
1 & 1 & 0 & 0 \\ 
0 & 1 & 1 & 0 \\ 
1 & 0 & 0 & 1 \\ 
0 & 1 & 0 & 1%
\end{array}%
\right) $ $A_{44} = \left( 
\begin{array}{cccc}
1 & 0 & 0 & 0 \\ 
1 & 0 & 1 & 0 \\ 
0 & 1 & 1 & 0 \\ 
1 & 0 & 0 & 1%
\end{array}%
\right) $ $A_{45} = \left( 
\begin{array}{cccc}
0 & 1 & 0 & 1 \\ 
1 & 0 & 1 & 0 \\ 
0 & 0 & 1 & 1 \\ 
0 & 1 & 0 & 0%
\end{array}%
\right) $ $A_{46} = \left( 
\begin{array}{ccc}
1 & 0 & 0 \\ 
0 & 0 & 1 \\ 
1 & 1 & 0 \\ 
0 & 1 & 1%
\end{array}%
\right) $ $A_{47} = \left( 
\begin{array}{ccc}
1 & 0 & 0 \\ 
1 & 1 & 0 \\ 
0 & 1 & 1%
\end{array}%
\right) $ $A_{48} = \left( 
\begin{array}{cc}
1 & 0 \\ 
1 & 1 \\ 
0 & 1%
\end{array}%
\right) $ $A_{49} = \left( 
\begin{array}{cc}
1 & 0 \\ 
1 & 1%
\end{array}%
\right) $

$A_{50} = \left( 
\begin{array}{cc}
1 & 0 \\ 
1 & 1%
\end{array}%
\right) $ $A_{51} = \left( 
\begin{array}{cc}
1 & 1 \\ 
1 & 0%
\end{array}%
\right) $ $A_{52} = \left( 
\begin{array}{c}
1 \\ 
1%
\end{array}%
\right) $ $A_{53} = \left( 
\begin{array}{c}
1%
\end{array}%
\right) $

$A_{54}=\left( 
\begin{array}{c}
1%
\end{array}%
\right) $ $A_{55}=\left( 
\begin{array}{c}
1%
\end{array}%
\right) $ $A_{56}=\left( 
\begin{array}{c}
1%
\end{array}%
\right) $ $A_{57}=\left( 
\begin{array}{c}
1%
\end{array}%
\right) .$

\textbf{The additive basis of} $H^{\ast }(E_{8}/E_{7})$ derived from the $A_{%
{k}}^{\prime }s$

%

\textbf{7.3. The ring structure on} $H^{even}(E_{8}/E_{7})$: See in the last
column of the table corresponding to $H^{even}$.

\bigskip

\textbf{7.4. The structure matrices }$M(\pi _{{m}})$\textbf{\ }are with
respect to the ordered monomial base $B(2m)$ (cf. [DZ$_{3}$, 3.3])

$B(30)=\{ {y_{15}}, {y_1}^{3}{y_{6}}^{2}, {y_1}^{5}{y_{10}}, {y_1}^{9}{y_{6}}%
, {y_1}^{15} \}$

$B(40)=\{ {y_{10}}^{2}, {y_1}^{2}{y_{6}}^{3}, {y_1}^{4}{y_{6}}{y_{10}}, {y_1}%
^{5}{y_{15}}, {y_1}^{8}{y_{6}}^{2}, {y_1}^{10}{y_{10}}, {y_1}^{14}{y_{6}}, {%
y_1}^{20} \}$

$B(48)=\{ {y_{6}}^{4}, {y_1}^{2}{y_{6}}^{2}{y_{10}}, {y_1}^{3}{y_{6}}{y_{15}}%
, {y_1}^{4}{y_{10}}^{2}, {y_1}^{6}{y_{6}}^{3}, {y_1}^{8}{y_{6}}{y_{10}}, {y_1%
}^{9}{y_{15}},$

$\hspace{2cm} {y_1}^{12}{y_{6}}^{2}, {y_1}^{14}{y_{10}}, {y_1}^{18}{y_{6}}, {%
y_1}^{24}\}$

$B(60)=\{ {y_{15}}^{2}, {y_{10}}^{3}, {y_{6}}^{5}, {y_1}^{2}{y_{6}}^{3}{%
y_{10}}, {y_1}^{3}{y_{6}}^{2}{y_{15}}, {y_1}^{4}{y_{6}}{y_{10}}^{2}, {y_1}%
^{5}{y_{10}}{y_{15}}, {y_1}^{6}{y_{6}}^{4}, $

$\hspace{2cm} {y_1}^{8}{y_{6}}^{2}{y_{10}}, {y_1}^{9}{y_{6}}{y_{15}}, {y_1}%
^{10}{y_{10}}^{2}, {y_1}^{12}{y_{6}}^{3}, {y_1}^{14}{y_{6}}{y_{10}}, {y_1}%
^{15}{y_{15}},$

$\hspace{2cm} {y_1}^{18}{y_{6}}^{2}, {y_1}^{20}{y_{10}}, {y_1}^{24}{y_{6}}, {%
y_1}^{30} \}$

\noindent of\textsl{\ }$\mathbb{Z}[y_{{1}},y_{6},y_{{10}},y_{{15}}]^{(2m)}$
for $m=15,20,24,30$ respectively, where $y_{{1}}=s_{{1,1}},y_{{6}%
}=s_{6,2},y_{{10}}=s_{10,1},y_{{15}}=s_{15,4}$.

$M( \pi_{15}) = \left( 
\begin{array}{cccc}
0 & 0 & 0 & 1 \\ 
36 & 26 & 35 & 32 \\ 
10 & 7 & 10 & 9 \\ 
63 & 45 & 62 & 56 \\ 
110 & 78 & 110 & 98 \\ 
&  &  & 
\end{array}%
\right) $

$M( \pi_{20}) = \left( 
\begin{array}{cccccc}
17 & 14 & 49 & 16 & 41 & 33 \\ 
384 & 324 & 1113 & 366 & 920 & 744 \\ 
106 & 89 & 309 & 101 & 257 & 207 \\ 
5 & 5 & 20 & 7 & 20 & 15 \\ 
670 & 564 & 1945 & 638 & 1611 & 1301 \\ 
185 & 155 & 540 & 176 & 450 & 362 \\ 
1169 & 982 & 3399 & 1112 & 2821 & 2275 \\ 
2040 & 1710 & 5940 & 1938 & 4940 & 3978 \\ 
&  &  &  &  & 
\end{array}%
\right) $

$M( \pi_{24}) = \left( 
\begin{array}{ccccccc}
7308 & 12639 & 17970 & 1422 & 6560 & 5406 & 9204 \\ 
2023 & 3512 & 4987 & 392 & 1822 & 1496 & 2560 \\ 
134 & 242 & 337 & 26 & 127 & 98 & 182 \\ 
560 & 976 & 1384 & 108 & 506 & 414 & 712 \\ 
12756 & 22092 & 31396 & 2478 & 11464 & 9435 & 16093 \\ 
3531 & 6139 & 8713 & 683 & 3184 & 2611 & 4476 \\ 
234 & 423 & 589 & 45 & 222 & 171 & 318 \\ 
22265 & 38616 & 54853 & 4318 & 20034 & 16467 & 28138 \\ 
6163 & 10731 & 15223 & 1190 & 5564 & 4557 & 7826 \\ 
38862 & 67500 & 95836 & 7524 & 35010 & 28740 & 49198 \\ 
67830 & 117990 & 167440 & 13110 & 61180 & 50160 & 86020 \\ 
&  &  &  &  &  & 
\end{array}%
\right) $

$M( \pi_{30}) = \left( 
\begin{array}{ccccccc}
135 & 275 & 295 & 329 & 201 & 75 & 15 \\ 
7911 & 16181 & 17399 & 19583 & 12117 & 4657 & 1123 \\ 
645195 & 1321140 & 1420270 & 1600126 & 991641 & 382055 & 92450 \\ 
179208 & 366799 & 394356 & 444129 & 275076 & 105884 & 25589 \\ 
12388 & 25139 & 27084 & 30331 & 18650 & 7114 & 1704 \\ 
49776 & 101838 & 109498 & 123272 & 76304 & 29344 & 7082 \\ 
3440 & 6980 & 7520 & 8419 & 5173 & 1970 & 470 \\ 
1127568 & 2308528 & 2481820 & 2795728 & 1732212 & 667154 & 161360 \\ 
313190 & 640937 & 689109 & 775979 & 480505 & 184895 & 44661 \\ 
21648 & 43928 & 47327 & 52995 & 32577 & 12419 & 2971 \\ 
86990 & 177950 & 191340 & 215380 & 133288 & 51240 & 12360 \\ 
1970580 & 4033869 & 4336803 & 4884675 & 3025848 & 1164996 & 281631 \\ 
547341 & 1119961 & 1204169 & 1355785 & 839349 & 322863 & 77947 \\ 
37830 & 76760 & 82700 & 92594 & 56904 & 21680 & 5180 \\ 
3443856 & 7048696 & 7578254 & 8534468 & 5285580 & 2034328 & 491542 \\ 
956550 & 1957000 & 2104200 & 2368818 & 1466178 & 563780 & 136040 \\ 
6018600 & 12316750 & 13242460 & 14911360 & 9232890 & 3552350 & 857900 \\ 
10518300 & 21522060 & 23140260 & 26053020 & 16128060 & 6203100 & 1497300 \\ 
&  &  &  &  &  & 
\end{array}%
\right) $

\textbf{7.5. The nullspaces of }$M(\pi _{{m}})$ for $m=15,20,24,30$ are
respectively

$N( \pi_{15}) = \left( 
\begin{array}{ccccc}
-2 & 10 & 16 & -10 & 1 \\ 
&  &  &  & 
\end{array}%
\right) $

$N( \pi_{20}) = \left( 
\begin{array}{cccccccc}
-30 & -100 & -180 & 18 & 90 & -24 & 0 & 1 \\ 
-3 & -10 & -18 & 2 & 8 & -4 & 1 & 0 \\ 
&  &  &  &  &  &  & 
\end{array}%
\right) $

$N( \pi_{24}) = \left( 
\begin{array}{ccccccccccc}
-180 & -1080 & -30 & -525 & 200 & 330 & 60 & 0 & 0 & 0 & 1 \\ 
-100 & -600 & -12 & -285 & 110 & 186 & 30 & 0 & 0 & 1 & 0 \\ 
-90 & -540 & -10 & -255 & 100 & 170 & 26 & 0 & 2 & 0 & 0 \\ 
-20 & -120 & -2 & -57 & 20 & 34 & 6 & 2 & 0 & 0 & 0 \\ 
&  &  &  &  &  &  &  &  &  & 
\end{array}%
\right) $

{\small \setlength{\arraycolsep}{0.1 mm} 
\begin{equation*}
N( \pi_{30}) = \left( 
\begin{array}{cccccccccccccccccc}
-9641868 & -6989880 & 3942456 & -24909120 & -10590840 & -172727040 & \
89697384 & 0 &  &  &  &  &  &  &  &  &  &  \\ 
-3024855 & -4831025 & 723120 & -17910240 & -2088060 & -65062950 & 28165110 \ 
& 0 &  &  &  &  &  &  &  &  &  &  \\ 
-694389 & -1773945 & 598104 & -2902180 & -817650 & -16256835 & 6431232 & 0 & 
&  &  &  &  &  &  &  &  &  \\ 
-765732 & -1490018 & -586714 & -11067120 & 141246 & -21238758 & 7264488 & 0
&  &  &  &  &  &  &  &  &  &  \\ 
919569 & -981465 & 206890 & -1456400 & -690090 & 2702760 & -8204562 & 0 &  & 
&  &  &  &  &  &  &  &  \\ 
-160808 & -592113 & 189928 & -1906990 & -158024 & -6447138 & 1670422 & 0 & 
&  &  &  &  &  &  &  &  &  \\ 
-148023 & -79164 & -958773 & -6254040 & 391392 & -4362987 & 1294896 & 0 &  & 
&  &  &  &  &  &  &  &  \\ 
-53081 & -56215 & 73266 & -295020 & -45720 & -1620600 & 376802 & 0 &  &  & 
&  &  &  &  &  &  &  \\ 
579307 & -196293 & 41378 & -291280 & -2114984 & 540552 & -4804058 & 0 &  & 
&  &  &  &  &  &  &  &  \\ 
-12879 & -649549 & 25784 & -2289480 & -6214 & -4456989 & 453438 & 0 &  &  & 
&  &  &  &  &  &  &  \\ 
-305733 & 2947385 & -4603644 & -16578720 & 922710 & 6632400 & 1028484 & 
4942415 &  &  &  &  &  &  &  &  &  &  \\ 
&  &  &  &  &  &  &  &  &  &  &  &  &  &  &  &  & 
\end{array}%
\right.
\end{equation*}
\begin{equation*}
\left. 
\begin{array}{cccccccccccccccccc}
0 & 0 & 0 & 0 & 0 & 0 & 0 & 0 & 0 & 988483 &  &  &  &  &  &  &  &  \\ 
0 & 0 & 0 & 0 & 0 & 0 & 0 & 0 & 988483 & 0 &  &  &  &  &  &  &  &  \\ 
0 & 0 & 0 & 0 & 0 & 0 & 0 & 988483 & 0 & 0 &  &  &  &  &  &  &  &  \\ 
0 & 0 & 0 & 0 & 0 & 0 & 988483 & 0 & 0 & 0 &  &  &  &  &  &  &  &  \\ 
0 & 0 & 0 & 0 & 0 & 988483 & 0 & 0 & 0 & 0 &  &  &  &  &  &  &  &  \\ 
0 & 0 & 0 & 0 & 988483 & 0 & 0 & 0 & 0 & 0 &  &  &  &  &  &  &  &  \\ 
0 & 0 & 0 & 988483 & 0 & 0 & 0 & 0 & 0 & 0 &  &  &  &  &  &  &  &  \\ 
0 & 0 & 988483 & 0 & 0 & 0 & 0 & 0 & 0 & 0 &  &  &  &  &  &  &  &  \\ 
0 & 1976966 & 0 & 0 & 0 & 0 & 0 & 0 & 0 & 0 &  &  &  &  &  &  &  &  \\ 
1976966 & 0 & 0 & 0 & 0 & 0 & 0 & 0 & 0 & 0 &  &  &  &  &  &  &  &  \\ 
0 & 0 & 0 & 0 & 0 & 0 & 0 & 0 & 0 & 0 &  &  &  &  &  &  &  &  \\ 
&  &  &  &  &  &  &  &  &  &  &  &  &  &  &  &  & 
\end{array}%
\right)
\end{equation*}
}

\begin{center}
\textbf{References}
\end{center}

\begin{enumerate}
\item[{[DZ$_{1}$]}] H. Duan and Xuezhi Zhao, Algorithm for multiplying
Schubert classes. Internat. J. Algebra Comput. 16(2006), 1197--1210.

\item[{[DZ$_{3}$]}] H. Duan and Xuezhi Zhao, The Chow rings of generalized
Grassmannians, \textsl{available on} arXiv: math.AG/0511332

\item[{[Hu]}] J. E. Humphreys, Introduction to Lie algebras and
representation theory, Graduated Texts in Math. 9, Springer-Verlag New York,
1972.
\end{enumerate}

\end{document}